\input ./style/arxiv-general.cfg
\documentclass[aap,MSNbibl,seceqn,dvips]{arximspdf}
\makeatletter
   \@ifpackageloaded{graphicx}{}{\usepackage{graphicx}}
\makeatother

\doi{10.1214/14-AAP1050} 
\volume{25}
\issue{4}
\pubyear{2015}
\firstpage{2339}
\lastpage{2381}
\docsubty{FLA}

\makeatletter
\newcommand{\rrvert}{\vert}
\newcommand{\llvert}{\vert}
\def\bold{\mathbf}
\newcommand{\eqref}[1]{(\ref{#1})}
\newtheorem{theorem} {Theorem}[section]
\newtheorem{lemma}[theorem]{Lemma}

\newtheorem{proposition}[theorem]{Proposition}
\newproclaim{remark}[theorem]{Remark}
\newproclaim{rem}{Remark}
\newcommand{\Dim}{\operatorname{Dim}}

\makeatother

\begin{document}
\begin{frontmatter}

\title{Limit shapes for growing extreme characters~of~$U(\infty)$}
\runtitle{Limit shapes for growing extreme characters of $U(\infty)$}

\begin{aug}
\begin{aug}
\author[A]{\fnms{Alexei}~\snm{Borodin}\ead[label=e1]{borodin@math.mit.edu}\thanksref{T1}},
\author[B]{\fnms{Alexey}~\snm{Bufetov}\corref{}\ead[label=e2]{alexey.bufetov@gmail.com}\thanksref{T2}}
\and
\author[B]{\fnms{Grigori}~\snm{Olshanski}\ead[label=e3]{olsh2007@gmail.com}\thanksref{T3}}
\thankstext{T1}{Supported in part by the NSF Grant DMS-10-56390.}
\thankstext{T2}{Supported in part by Simons Foundation-IUM
scholarship, by Moebius Foundation for Young Scientists, by ``Dynasty''
Foundation and by the RFBR Grant 13-01-12449.}
\thankstext{T3}{Supported in part by a grant from Simons foundation
(Simons-IUM Fellowship), by the project SFB 701 of Bielefeld University
and by the RFBR Grant 13-01-12449.}
\affiliation{Massachusetts Institute of Technology and Institute for Information Transmission Problems,
Institute for Information Transmission Problems and Higher School of
Economics, and
Institute for Information Transmission Problems and Higher School of Economics}
\runauthor{A. Borodin, A.~Bufetov and G. Olshanski}
\address[A]{A. Borodin\\
Department of Mathematics\\
Massachusetts Institute of Technology\\
Cambridge, MA\\
USA\\
and\\
Institute for Information\\
\quad Transmission Problems\\
Moscow\\
Russia\\
\printead{e1}}
\address[B]{A. Bufetov\\
G. Olshanski\\
Institute for Information\\
\quad Transmission Problems\\
Moscow\\
Russia\\
and\\
Higher School of Economics\\
Moscow\\
Russia\\
\printead{e2}\\
\phantom{E-mail:\ }\printead*{e3}}
\end{aug}

\end{aug}

\received{\smonth{12} \syear{2013}}
\revised{\smonth{6} \syear{2014}}

%
\begin{abstract}
We prove the existence of a limit shape and give its explicit
description for certain probability
distribution on signatures (or highest weights for unitary groups). The
distributions have representation
theoretic origin---they encode decomposition on irreducible characters
of the restrictions
of certain extreme characters of the infinite-dimensional unitary
group $U(\infty)$ to growing finite-dimensional unitary subgroups
$U(N)$. The
characters of $U(\infty)$ are allowed to depend on $N$.
In a special case, this describes the hydrodynamic behavior for a
family of
random growth models in $(2+1)$-dimensions with varied initial conditions.
\end{abstract}

%


\end{frontmatter}

\section{Introduction}\label{sec1}

Decomposing the restriction of an irreducible representation
of a group to its subgroup onto irreducible components is one of the basic
problems of the representation theory. Under special circumstances, as the
group and the subgroup become large, such decomposition may be
subject to a Law of Large Numbers type concentration phenomenon---the
bulk of
the decomposition consists of representations that are in some sense
close to
each other. This paper is devoted to studying one of such situations.

Historically, the first example of this concentration phenomenon was discovered
by Vershik--Kerov \cite{VK-77} and Logan--Shepp \cite{LS}. One way to phrase
their result is to consider the infinite bisymmetric group $G=S(\infty
)\times
S(\infty)$, where $S(\infty)$ is the group of finite permutations of
$\mathbb
N:=\{1,2,\ldots\}$, and the growing subgroups being finite bisymmetric groups
$G(n)=S(n)\times S(n)$, where $S(n)$ consists of permutations of a
subset of
$\Bbb N$ with $n\ge1$ elements. Take the biregular representation of
$G$ in
$\ell^2(S(\infty))$ with $G$ acting by left and right shifts. It is
well-known
that it is irreducible (as for any countable group with infinite nontrivial
conjugacy classes). Its restriction to $G(n)$ decomposes on isotypical
components corresponding to irreducible representations of $S(n)$, or to
partitions of $n$ (equivalently, Young diagrams with $n$ boxes). The
corresponding spectral measure is the celebrated \textit{Plancherel distribution}
on Young diagrams with $n$ boxes that assigns to $\lambda$ the weight
equal to
the square of the number of standard Young tableaux of shape $\lambda$ divided
by $n!$.

The theorem of Vershik--Kerov--Logan--Shepp (see Kerov \cite{K-CLT} and
Ivanov--Olshanski \cite{IvaOls} for a different proof that is closer to
the present
work) says that if we shrink the random Young diagram $\lambda$ by the factor
of $\sqrt n$ in both directions (so that its area is now 1), then as
$n\to
\infty$, the boundary of $\lambda$ converges, in probability and in a suitable
topology, to an explicit smooth curve usually referred to as the \textit{limit
shape}.

Vershik--Kerov in \cite{VK-81} also considered the case of other
(unitary spherical) irreducible representations of $G$ and their
restrictions to $G(n)$, showing that while the law of large numbers is still
there, it takes a drastically different form---one needs to normalize
the row
and column lengths of the random Young diagram $\lambda$ by $n$ to see the
almost sure convergence to a point configuration (not a smooth curve) that
essentially encodes the original representation of $G$.

\begin{figure}[b]

\includegraphics{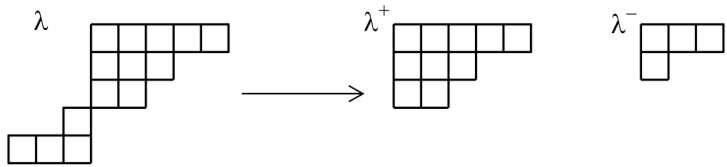}

\caption{Signature $\lambda= (5,3,2,-1,-3)$ and Young diagrams $\lambda
^+=(5,3,2)$ and $\lambda^-=(3,1)$.}
\label{2Yd}
\end{figure}

In the present paper, we are dealing not with the symmetric groups
$S(n)$ but
with the compact unitary groups $U(N)$. The irreducible characters of $U(N)$
are parameterized by $N$-tuples $\lambda=(\lambda_1\ge\lambda_2\ge\cdots
\ge
\lambda_N)\in\Bbb Z^N$, which are called \emph{signatures} of length
$N$. Note
that every such $\lambda$ can be viewed as a couple $(\lambda^+,\lambda
^-)$ of
Young diagrams (their row-lengths are, resp., the positive and the minus
negative coordinates in $\lambda$; see an example in Figure~\ref{2Yd}).
These two Young diagrams represent the
\emph{shape} of the signature.

Let us take $G(N)=U(N)\times U(N)$ and define $G$ as the union of the growing
groups $G(N)$. In other words, $G=U(\infty)\times U(\infty)$, where
$U(\infty)$
is the group of unitary matrices of format $\mathbb N\times\Bbb N$ with
finitely many entries $U_{ij}$ distinct from $\delta_{ij}$. The
restriction of
a (unitary spherical) irreducible representation of $G$ to $G(N)$
decomposes on
isotypical components parameterized by the signatures of length~$N$. It is
easiest to encode this decomposition via \emph{characters}---central normalized positive-definite functions on $U(\infty)$ that
are in
one-to-one correspondence with the spherical unitary representations; see
Olshanski \cite{Olsh1,Olsh2}.

If $\chi\dvtx U(\infty)\to\mathbb C$ is a character of $U(\infty)$, then
\[
\chi\bigl(\operatorname{diag}(z_1,\ldots,z_N,1,1,
\ldots)\bigr)=\sum_{
\lambda=(\lambda_1\ge\cdots\ge\lambda_N)\in\mathbb Z^N} {M_N^\chi}(
\lambda)\frac{s_\lambda(z_1,\ldots,z_N)}{s_\lambda(1,\ldots
,1)},
\]
where $s_\lambda$'s are the rational Schur functions [conventional irreducible
characters for $U(N)$], and $M_N^\chi$ is the spectral measure of the
decomposition, which is a probability distribution on the set of all signatures
of length $N$.

Irreducible (spherical unitary) representations of $G$ correspond to the
extreme points of the convex set of characters of $U(\infty)$, often referred
to as its \emph{extreme characters}. The classification of the extreme
characters is known as the Edrei--Voiculescu theorem (see Voiculescu
\cite{Vo},
Edrei \cite{Edr}, Vershik--Kerov \cite{VK}, Okounkov--Olshanski \cite{OkoOls98},
Borodin--Olshanski \cite{BorOls}). They can be parameterized by the set
\[
\Omega=\bigl(\alpha^+, \alpha^-, \beta^+,\beta^-, \delta^+, \delta^-\bigr)\in
\bigl(\Bbb R_+^\infty\bigr)^4\times(\Bbb R_+)^2,
\]
where
\begin{eqnarray*}
\alpha^{\pm} &=& \alpha_1^{\pm} \ge
\alpha_2^{\pm} \ge\cdots\ge0,\qquad \beta^{\pm} =
\beta_1^{\pm} \ge\beta_2^{\pm} \ge
\cdots\ge0,
\\
\delta^{\pm} &\ge&0,\qquad \sum_{i=1}^{\infty}
\bigl(\alpha_i^{\pm}+\beta _i^{\pm}
\bigr) \le \delta^{\pm}, \qquad \beta_1^+ + \beta_1^-
\le1.
\end{eqnarray*}
Instead of $\delta^{\pm}$, we will use parameters $\gamma^{\pm}\ge0$
defined by
$
\gamma^{\pm}:= \delta^{\pm} - \sum_{i=1}^{\infty} (\alpha_i^{\pm} +
\beta_i^{\pm})$.
Each $\omega\in\Omega$ defines a function $\Phi^{\omega} \dvtx \{ u \in
\Bbb C \dvtx |u|=1 \} \to\Bbb C$ by
%
\begin{eqnarray}
\label{eq4} \Phi^{\omega} (u) &= &\exp\bigl( \gamma^+ (u-1) +
\gamma^{-} \bigl(u^{-1}-1\bigr) \bigr)
\nonumber
\\[-8pt]
\\[-8pt]
\nonumber
&&{}\times
\prod
_{i=1}^{\infty} \frac{(1+\beta_i^+ (u-1))}{(1-\alpha_i^+ (u-1))} \frac{(1+\beta_i^- (u^{-1}-1))}{(1- \alpha_i^- (u^{-1}-1))} ,
\end{eqnarray}
which we call the \emph{Voiculescu function} with parameter $\omega$. The
corresponding extreme character has the form
%
\begin{equation}
\label{eq5} \chi^{\omega} (U) := \prod_{u \in\operatorname{Spectrum}(U)}
\Phi ^{\omega} (u),\qquad U \in U(\infty),
\end{equation}
where the product is over the eigenvalues of $U$ [this product is essentially
finite, because $\Phi^{\omega}(1)=1$ and only finitely many of $u$'s are
distinct from $1$].

We are thus interested in the limit shape phenomenon for the probability
measures of the form $M_N^{\chi^\omega}$ as $N\to\infty$.

Let $\lambda=\lambda(N)$ be the random signature with distribution
$M_N^{\chi^\omega}$, and let $\lambda^\pm$ be the corresponding Young diagrams.
The row and column lengths of $\lambda^\pm$ (see Figure~\ref{2Yd})
divided by
$N$ almost surely converge, as $N\to\infty$, to the values of the
$\alpha^\pm$
and $\beta^\pm$ coordinates (somewhat similarly to the case of $S(\infty
)$, cf.
Vershik--Kerov~\cite{VK-81}). If all those coordinates are zero but
$\gamma^\pm$
are not, then scaling by $\sqrt{N}$ leads to concentration of
$M_N^{\chi^\omega}$ around two copies of the
Vershik--Kerov--Logan--Shepp limit
shape; see Borodin--Kuan \cite{BorKuan}.\setcounter{footnote}{3}\footnote{There is no proof of the
measure concentration there, but there is substantial evidence that it holds.}
The latter work also noted a hypothetical limit shape formation as
$\gamma^\pm$
\emph{grow linearly in $N$} as $N\to\infty$ (as opposed to being
independent of
$N$), and suggested a formula for the limit shape. In the case when only
$\gamma^+$ is nonzero, the concentration around the limit shape was proved
earlier by Biane \cite{Bia01}.

In the present work, we prove that the limit shape phenomenon takes
place in a
much more general setting.

Let us state our main result.

Consider a sequence of points $\omega(N)\in\Omega$, $N\ge1$, and
assume that
there exists an analytic function $P(z)$ defined in a neighborhood of the
origin such that
%
\begin{equation}
\label{11} \lim_{N\to\infty} \frac{1}N \bigl(\log
\Phi^{\omega(N)}(z+1) \bigr)=P(z)
\end{equation}
uniformly in a (possibly smaller) neighborhood of $z=0$ (see the
beginning of
Section~\ref{sec3} below for simple sufficient conditions for the above
convergence to
hold).

\begin{theorem}\label{thm1}
Let us fix an arbitrary sequence $\{\omega(N)\}_{N\ge1}$ of elements in
$\Omega$ satisfying the limit relation \eqref{11}. For every $N$, let
$\lambda(N)$ denote the random signature distributed according to
$M_N^{\chi^{\omega(N)}}$ and let $\lambda^\pm(N)$ be the corresponding Young
diagrams.

Let us shrink the diagrams $\lambda^\pm(N)$ by the factor of $N$ in both
directions. Then the resulting random shapes converge, as $N\to\infty$, to
certain nonrandom shapes, which in principle can be obtained from the function
$P(z)$.
\end{theorem}

Here is another (and more precise) formulation of the result.

Denote by $\delta(x)$ the Dirac measure at a point $x\in\mathbb R$. To every
signature $\lambda=(\lambda_1\ge\cdots\ge\lambda_N)$, we assign an atomic
probability measure on $\Bbb R$:
%
\begin{equation}
\label{eq6} \mu_\lambda:=\frac{1}N \sum
_{i=1}^N \delta \biggl(\frac{\lambda_i-i+1/2}{N} \biggr).
\end{equation}
This measure encodes the (scaled) shape of $\lambda$ (see Section~\ref{graph-sign} for more detail).

\begin{theorem}\label{thm2}
Let $\{\omega(N)\}_{N\ge1}$and $\lambda(N)$ be as above, and let
$\mu_{\lambda(N)}$ be the random atomic measure on $\mathbb R$
corresponding to
$\lambda(N)$.

There exists a probability measure $\sigma$ with compact support on
$\mathbb R$
such that
\[
\lim_{N\to\infty} \mu_{\lambda(N)} = \sigma\qquad \mbox{(weak
convergence in probability)}.
\]
The measure $\sigma$ is uniquely determined by its moments $(1,m_1,m_2,
\ldots)$, which in turn are found from the fact that the two formal
series in
$z$,
\[
\exp\bigl(z+m_1z+m_2z^3+\cdots\bigr)-1\quad
\mbox{and}\quad \frac{z}{1+z(1+z)P'(z)}
\]
are mutually inverse with respect to composition.
\end{theorem}

See Theorem~\ref{mainTh} below. The fact that Theorem~\ref{thm2} implies
Theorem~\ref{thm1} is explained in Proposition~\ref{prop2}. That proposition
also shows that the limit measure $\sigma$ always has a density with
respect to
Lebesgue measure.

Concrete example of sequences $\{\omega(N)\}_{N\ge1}$ and
corresponding limit
shapes can be found in the \hyperref[app]{Appendix} below.

The density of $\sigma$ can be guessed using the determinantal
structure of
suitably defined correlation functions of measures $M_N^\omega$ found by
Borodin--Kuan \cite{BorKuan}, and a steepest descent analysis of the double
contour integral representation of the correlation kernel. We outline this
route in Section~\ref{DetPr} below. Note, however, that \emph{proving} the
concentration of
measure phenomenon is a different task, and correlation functions are
not well
suited for it. In this work, we employ a different approach.

Our result also has a probabilistic interpretation. Measures of the form
$M_N^\omega$ with $\omega$ having finitely many nonzero $\alpha^\pm$ and
$\beta^\pm$ parameters can be obtained via a Markov growth process in
$(2+1)$-dimensions; see Borodin--Ferrari \cite{BorFer}. Our main result then
establishes the law of large numbers for a growing two-dimensional random
interface. The growth process is \emph{local}, and one can expect that the
limit shape should be evolving in time according to a first-order PDE. Our
result confirms that for a broad class of initial conditions; see
Section~\ref{sec3.3} for
details.\footnote{In the case when the only nonzero paramater is $\gamma^+$,
the corresponding PDE was found in \cite{BorFer}.}

If we have two sequences of extreme characters that lead to limit
shapes, we
can also consider the sequence whose members are products of those of
the two
original sequences (the set of extreme characters is closed under
multiplication). The new sequence will also have a limit shape, and we thus
obtain an operation on limiting measures $\sigma$. We call it
``quantized free
convolution''; it is a relative of the free convolution in free
probability, and
it degenerates to it; see Section~\ref{sec3.4} below. Bufetov--Gorin \cite{BufGor} show how this
operation naturally arises through tensoring large irreducible representations
of growing (but finite-dimensional) unitary groups and further decomposing
them on irreducibles.

The particular examples of characters of $U(\infty)$ are the one-sided
Plancherel character (the only nonzero parameter is $\gamma^+$) and the
two-sided Plancherel chararcter ($\gamma^+$ and $\gamma^-$ are
nonzero). The probability measures arising from these characters were
considered, for example, by Biane \cite{Bia01}, Borodin--Bufetov \cite{BorBuf}, Borodin--Kuan \cite{BorKuan}. However, we want to emphasize
that the conditions of Theorem~\ref{thm2} are much more general because
they allow to manipulate not only $\gamma^+$, $\gamma^-$, but all $4
\cdot\infty+2$ parameters of extremal characters. Theorem~\ref{thm2}
gives the same answer that was proved earlier by Biane \cite{Bia01} in
the case of the one-sided Plancherel character and conjectured by
Borodin--Kuan \cite{BorKuan} in the case of the two-sided Plancherel character.

Having proved a law of large numbers, it is natural to ask about the central
limit theorem. In the case of linearly growing parameter $\gamma^+$ and all
other parameters being zero, it was shown in Borodin--Ferrari \cite{BorFer} and
Borodin--Bufetov \cite{BorBuf} that the fluctuations around the limit
shape are
described by the \emph{two-dimensional Gaussian Free Field.} It is plausible
that a similar description of fluctuations should exist under the
(substantially more general) assumption of our theorem above.

Our proof is based on the method of moments. It bears a certain
similarity with
the work of Ivanov--Olshanski \cite{IvaOls} for the Plancherel measures
on symmetric
groups and the work of Borodin--Bufetov \cite{BorBuf} for the nonzero
$\gamma^+$ case,
but it is of course more involved because of the many parameters
present. The
key ingredients are provided by certain graph enumeration arguments, as we
explain in Section~\ref{proofTh}.

\section{Preliminaries}\label{sec2}

\subsection{The infinite-dimensional unitary group and its characters}
\label{infinite group}

Let\break $U(N)= \{ [u_{ij}]_{i,j=1}^N  \}$ be the group of $N
\times N$ unitary matrices.
Consider the tower of embedded unitary groups\vspace*{2pt}
\[
U(1) \subset U(2) \subset\cdots\subset U(N) \subset U(N+1) \subset \cdots,
\]
where the embedding $U(N) \subset U(N+1)$ is defined by
$u_{i,N+1}=u_{N+1,i}=0$, $1 \le i \le k$, $u_{N+1,N+1}=1$. \emph{The
infinite-dimensional unitary group} is the union of these groups:\vspace*{2pt}
\[
U(\infty) = \bigcup_{N=1}^{\infty} U(N).
\]

Define a \emph{character} of the group $U(\infty)$ as a function $\chi\dvtx
U(\infty) \to\mathbb C$ that satisfies the following conditions:
\begin{longlist}[(1)]
\item[(1)] $\chi(e)=1$, where $e$ is the identity element of $U(\infty)$
(normalization);

\item[(2)] $\chi(g h g^{-1}) = \chi(h)$, where $g,h$ are any elements of
$U(\infty)$ (centrality);

\item[(3)] $[\chi( g_i g_j^{-1})]_{i,j=1}^n$ is an Hermitian and
positive-definite matrix for
any $n\ge1$ and $g_1, \ldots, g_n \in U(\infty)$ (positive-definiteness);

\item[(4)] the restriction of $\chi$ to $U(N)$ is a continuous function for any
$N\ge1$
(continuity).
\end{longlist}

The set of characters of $U(\infty)$ is obviously convex. The extreme
points of
this set are called the \emph{extreme} characters; they replace irreducible
characters in this setting. The classification of the extreme
characters was described in the \hyperref[sec1]{Introduction}; see formulas \eqref{eq4}
and \eqref{eq5}
above.

\subsection{The Gelfand--Tsetlin graph and coherent systems of measures}
\label{gelfand-tsetlin}

A \emph{signature} (also called \emph{highest weight}) of length $N$ is
a sequence of $N$ weakly decreasing integers\vspace*{2pt}
\[
\lambda= ( \lambda_1 \ge\lambda_2 \ge\cdots\ge
\lambda_N ),\qquad \lambda_i \in\mathbb Z, 1\le i\le N.
\]

It is well known that the irreducible (complex) representations of
$U(N)$ can
be parameterized by signatures of length $N$ (see, e.g., \cite{Wey39,Zhe}). Let $\Dim_N (\lambda)$ be the dimension of the representation
corresponding to $\lambda$. By $\chi^{\lambda}$, we denote the \emph{normalized}
character of this representation, that is, the conventional character divided
by $\Dim_N (\lambda)$.

Let $\mathbb{GT}_N$ denote the set of all signatures of length $N$. (Here,
letters $\mathbb{GT}$ stand for ``Gelfand--Tsetlin.'') We say that
$\lambda\in
\mathbb{GT}_N$ and $\mu\in\mathbb{GT}_{N-1}$ \emph{interlace},
notation $\mu
\prec\lambda$, if $\lambda_i \ge\mu_i \ge\lambda_{i+1}$ for any $1
\le i \le
N-1$. We also define $\mathbb{GT}_0$ as a singleton consisting of an element
that we denote as $\varnothing$. We assume that $\varnothing\prec
\lambda$ for
any $\lambda\in\mathbb{GT}_1$.

The \emph{Gelfand--Tsetlin graph} $\mathbb{GT}$ is defined by
specifying its set of vertices as
$\bigcup_{N=0}^{\infty} \mathbb{GT}_N $ and putting an edge between any
two signatures $\lambda$ and $\mu$ such that
either $\lambda\prec\mu$ or $\mu\prec\lambda$.
A \emph{path} between signatures $\kappa\in\mathbb{GT}_K$ and $\nu
\in\mathbb{GT}_N$, $K<N$, is a sequence
\[
\kappa= \lambda^{(K)} \prec\lambda^{(K+1)} \prec\cdots\prec
\lambda ^{(N)} = \nu,\qquad \lambda^{(i)} \in\mathbb{GT}_i,
K\le i\le N.
\]

It is well known that $\Dim_N (\nu)$ is equal to the number of paths
between $\varnothing$ and $\nu\in\mathbb{GT}_N$. An \emph{infinite
path} is a sequence
\[
\varnothing\prec\lambda^{(1)} \prec\lambda^{(2)} \prec\cdots
\prec \lambda^{(k)} \prec\lambda^{(k+1)} \prec\cdots.
\]

We denote by $\mathcal P$ the set of all infinite paths. It is a
topological space with
the topology induced from the product topology on the ambient product
of discrete sets
$\prod_{N\ge0}\mathbb{GT}_N$. Let us equip $\mathcal P$ with the Borel
$\sigma$-algebra.

For $N=0,1,2,\ldots,$ let $ M_N$ be a probability measure on $\mathbb
{GT}_N$. We say that
$\{ M_N \}_{N=0}^{\infty}$ \emph{is a coherent system of measures} if
for any $N\ge0$
and $\lambda\in\mathbb{GT}_N$,
\[
M_{N} (\lambda) = \sum_{\nu\dvtx  \lambda\prec\nu}
M_{N+1} (\nu) \frac
{\Dim_{N} (\lambda)}{\Dim_{N+1} (\nu)}.
\]

Given a coherent system of measures $\{ M_N \}_{N=1}^{\infty}$, define the
weight of a cylindric set of $\mathcal P$ consisting of all paths with
prescribed members up to $\mathbb{GT}_N$ by
%
\begin{equation}
\label{mera-puti} P \bigl( \lambda^{(1)}, \lambda^{(2)}, \ldots,
\lambda^{(N)} \bigr) = \frac{ M_N
(\lambda^{(N)})}{\Dim_N (\lambda^{(N)} )}.
\end{equation}
Note that this weight depends on $\lambda^{(N)}$ only (and does not
depend on $\lambda^{(1)},\break \lambda^{(2)}, \ldots,\lambda
^{(N-1)}$). The coherency property implies that
these weights are consistent, and they correctly define a Borel
probability measure on $\mathcal P$.

Now let $\chi$ be an arbitrary character of $U(\infty)$ and $\chi_N$
denote its
restriction to the subgroup $U(N)$. The function $\chi_N$ can be
expanded into
a series in $\chi^{\lambda}$'s,
%
\begin{equation}
\label{char} \chi_N = \sum_{\lambda\in\mathbb{GT}_N}
M_N (\lambda) \chi^{\lambda}.
\end{equation}
It is readily seen that the coefficients $M_N (\lambda)$ form a
coherent system
of measures on $\mathbb{GT}$. Conversely, for any coherent system of measures
on $\mathbb{GT}$ one can construct a character of $U(\infty)$ using the above
formula.

Note also that if $\chi_N$ is smooth, then the coefficients of the expansion
\eqref{char} rapidly decay as $\lambda$ goes to infinity, so that any
polynomial function in variables $\lambda_1,\ldots,\lambda_N$ is
summable on
$\mathbb{GT}_N$ with respect to measure $M_N$.

\subsection{The algebra of shifted symmetric functions}
\label{shifted-symmetric-functions}

In this subsection, we review some facts about the algebra of shifted symmetric
functions; see \cite{OkoOls,KerOls,IvaOls}.

Let $\operatorname{Sym}^* (N)$ be the algebra of polynomials in $N$ variables
$x_1, \ldots,x_N$, that are symmetric in shifted variables
\[
y_i:=x_i-i+\tfrac12,\qquad  i=1,2, \ldots, N.
\]

The standard filtration of $\operatorname{Sym}^*(N)$ is defined by the degree
of a polynomial. Define a map $\operatorname{Sym}^*(N) \to
\operatorname{Sym}^*(N-1)$ as specializing $x_N=0$. The \emph{algebra of
shifted symmetric functions} $\operatorname{Sym}^*$ is the projective
limit of
the algebras $\operatorname{Sym}^*(N)$ with respect to these maps.
Here, the
limit is taken in the category of filtered algebras meaning that the degree
does not grow.

The algebra $\operatorname{Sym}^*$ can be identified with the
subalgebra in
$\mathbb R[[x_1,x_2,\ldots]]$ generated by the algebraically independent system
$\{ \mathbf p_k \}_{k=1}^{\infty}$, where
\[
\mathbf p_k (x_1, x_2, \ldots) := \sum
_{i=1}^{\infty} \biggl( \biggl( x_i-i+
\frac12 \biggr)^k - \biggl(-i+\frac12 \biggr)^k \biggr),\qquad k=1, 2, \ldots.
\]

Let $\mathbb Y_n$ denote the set of partitions (or Young diagrams) $\nu
= (\nu_1
\ge\nu_2 \ge\cdots\ge0)$ with $|\nu| := \sum_{i\ge1} \nu_i = n$. Let
$\rho,
\nu\in\mathbb Y:=\mathbb Y_0 \cup\mathbb Y_1 \cup\mathbb Y_2 \cup
\cdots,$
and let $r=|\rho|$, \mbox{$n=|\nu|$}. For $r=n$, denote by $\psi_{\rho}^{\nu}$ the
value of the irreducible character of the symmetric group $S(n)$ corresponding
to $\nu$ on the conjugacy class indexed by $\rho$ (see, e.g., \cite{Mac,Sag} for details on symmetric groups). For $r<n$, denote by
$\psi_{\rho}^{\nu}$ the value of the same character on the conjugacy class
indexed by $\rho\cup1^{n-r} = (\rho, 1,1, \ldots,1) \in\mathbb Y_n$. Define
$p_{\rho}^{\#} \dvtx \mathbb Y \to\mathbb R$ by
\[
p_{\rho}^{\#} (\nu) = %
\cases{\displaystyle n(n-1) \cdots(n-r+1)
\frac{ \psi^{\nu}_{\rho}}{ \dim\nu}, &\quad $n \ge r$; \vspace*{2pt}
\cr
0, &\quad  $n <r$.} %
\]

Note that elements of $\operatorname{Sym}^*$ are well-defined functions
on the
set of all infinite sequences with finitely many nonzero terms. It
turns out
that there is a unique element $\mathbf p_{\rho}^{\#} \in\operatorname{Sym}^*$
such that $\mathbf p_{\rho}^{\#} (\nu)=p_{\rho}^{\#} (\nu)$ for all
$\nu\in
\mathbb Y$. It is known that the set $\{ \mathbf p_{\rho}^{\#} \}_{\rho
\in
\mathbb Y}$ is a linear basis in $\operatorname{Sym}^*$. When $\rho$ consists
of a single row, $\rho= (k)$, we denote the element $\mathbf p_{\rho
}^{\#}$ by
$\mathbf p_k^{\#}$. It is also known that the set $\{ \mathbf p_k^{\#}
\}_{k=1}^{\infty}$ is an algebraically independent system of generators of
$\operatorname{Sym}^*$. See \cite{IvaOls} for details.

The \emph{weight} of $\mathbf p_{\rho}^{\#}$ is defined by
\[
\operatorname{wt}\bigl( \mathbf p_{\rho}^{\#}\bigr) = |\rho|+l(\rho),
\]
where $l(\rho)$ denotes the number of nonzero coordinates in $\rho$. We extend
this definition to arbitrary elements $f\in\operatorname{Sym}^*$ in a natural
way, namely, we expand $f$ in the basis $\{ \mathbf p^\#_\rho\}$ and
define the
\emph{weight} $\operatorname{wt}(f)$ as the maximal weight of those basis elements
that enter
the expansion of $f$ with nonzero coefficients. It turns out (see~\cite{IvaOls}) that $\operatorname{wt}(\cdot)$ is a filtration on $\operatorname
{Sym}^*$. It
is called the \emph{weight filtration}.

We will need the following formula (see \cite{IvaOls}, Proposition~3.7):
%
\begin{equation}\qquad
\label{change} \mathbf p_k = \frac{1}{k+1} \bigl[
u^{k+1}\bigr] \bigl\{ \bigl(1 + \mathbf p_1^{\#}
u^2 + \mathbf p_2^{\#} u^3 + \cdots
\bigr)^{k+1} \bigr\} + \mbox{lower weight terms},
\end{equation}
where ``lower weight terms'' denotes terms with weight $\le k$, and
$[u^k]\{ A(u) \}$ stands for the coefficient of $u^k$ in a formal power
series $A(u)$.

\subsection{An algebra of functions on (random) signatures}
\label{functions-on-signatures}

In this section, we define an algebra of functions on the probability
space $(\mathbb{GT}_N, M_N)$ and state some properties of these functions.

For any $N \ge1$, define functions $p_k^{(N)}\dvtx \mathbb{GT}_N \to
\mathbb R$ by
%
\begin{eqnarray}
\label{functions-pk} p_k^{(N)} (\lambda) = \sum
_{i=1}^{N} \biggl( \biggl(\lambda_i - i +
\frac12 \biggr)^k - \Biggl(-i + \frac12 \Biggr)^k \biggr),
\nonumber
\\[-8pt]
\\[-8pt]
\eqntext{\lambda \in
\mathbb{GT}_N, 1 \le k \le N.}
\end{eqnarray}

Let $\mathbb A(N)$ be the algebra generated by $\{ p_k^{(N)} \}
_{k=1}^{N}$. It
is easy to see that for a fixed $N \ge1$, the functions $\{ p_k^{(N)}
\}_{k=1}^{N}$ are algebraically independent; therefore, they form a
system of
algebraically independent generators of $\mathbb A(N)$. Clearly, the algebras
$\mathbb A(N)$ and $\operatorname{Sym}^*(N)$ are naturally isomorphic.

Consider the map $\mathrm{pr}_N\dvtx \operatorname{Sym}^* \to\mathbb A(N)$ such
that $\mathrm{pr}_N
(\mathbf p_k) = p_{k}^{(N)}$. Denote by $p_{\rho}^{\# (N)}$ the
function $\mathrm{pr}_N
(\mathbf p_{\rho}^{\#})$.

Let $\chi$ be a character of $U(\infty)$, $\chi_N$ its restriction to $U(N)$,
and $M_N$ the corresponding probability measure on $\mathbb{GT}_N$, where
$N=1,2,\ldots .$ We consider the pair $(\mathbb{GT}_N, M_N)$ as a probability
space. Then the functions from $\mathbb A(N)$ turn into random
variables. Let
$\mathbf E_N$ be the expectation on this probability space. Note that
for any
$f \in\operatorname{Sym}^*$ we can consider the random variable $\mathrm{pr}_N (f)$.

With some ambiguity that should not lead to any confusion, we omit the index
$N$ in the notation of $p_{k}^{(N)}$ and $p_{\rho}^{\# (N)}$.

The complexification of $U(N)$ is the group $\operatorname
{GL}(N,\mathbb C)$,
which is an open subset of $\operatorname{Mat}(N,\mathbb C)$, the space of
$N\times N$ complex matrices. Let $x_{ij}$ be the natural coordinates in
$\operatorname{Mat}(N,\mathbb C)$ (where $1\le i,j\le N$) and $\partial_{ij}$
be the abbreviation for the (holomorphic) partial derivative operator
$\partial/\partial x_{ij}$. Note that any analytic function on the real
manifold $U(N)$ can be extended to a holomorphic function in a
neighborhood of
the identity matrix in $\operatorname{Mat}(N,\mathbb C)$.

\begin{proposition}\label{prop1} Assume that $\chi$ is such that for
every $N=1,2,\ldots,$ the function $\chi_N$ is analytic and so admits a
holomorphic extension to a neighborhood of $1$ in
$\operatorname{Mat}(N,\mathbb C)$. Then the following formula holds:
%
\begin{equation}
\label{difOp1} \qquad\mathbf E_N \bigl(p_{\rho}^{\#}
\bigr) = \sum_{1 \le i_1, \ldots, i_{|\rho|} \le N} \partial_{i_1 s(i_1)}
\,\partial_{i_2 s(i_2)} \cdots\partial_{i_{|\rho|}
i_{s(|\rho|)}} \chi_N (1+X)
\Big|_{X=0},
\end{equation}
where $s \in S(|\rho|)$ is an arbitrary permutation with cycle
structure $\rho$,
and $X=[x_{ij}]$ is a matrix from $\operatorname{Mat}(N,\mathbb C)$
close to
$0$.
\end{proposition}

Before proceeding to the proof, let us note that we will apply this
result only
to the extreme characters, and all such characters satisfy the
hypothesis of
the proposition, because every Voiculescu function is analytic.
However, there
exist nonextreme characters $\chi$ for which the functions $\chi_N$ are not
analytic and even not smooth.

\begin{pf*}{Proof of Proposition \ref{prop1}}
Because $\chi_N$ is analytic, all functions from $\mathbb A(N)$ are summable
with respect to $M_N$, so that the corresponding random variables have finite
expectation. Thus, the left-hand side of \eqref{difOp1} is well defined.

The key fact we need is Theorem~2 in Kerov--Olshanski \cite{KerOls}
(see also
Okounkov--Olshanski \cite{OkoOls}, Section~15). Here is its statement. Consider
the differential operator
\[
D_\rho=\sum_{\alpha_1,\ldots,\alpha_k, i_1,\ldots,i_k=1}^N
x_{\alpha
_1i_1},\ldots, x_{\alpha_ki_k} \partial_{\alpha_1i_{s(1)}}\cdots
\partial_{\alpha_ki_{s(k)}}
\]
on $\operatorname{Mat}(N,\mathbb C)$. Its restriction to the group
$\operatorname{GL}(N,\mathbb C)$ is invariant with respect to left and right
shifts, and one has
\[
D_\rho\chi^\lambda=p^{\#}_{\rho} (\lambda)
\chi^\lambda
\]
for any $\lambda\in\mathbb{GT}_N$. Let us recall that $\chi^\lambda$ denotes
the normalized irreducible character of $U(N)$ indexed by $\lambda$, so that
$\chi^\lambda(1)=1$. Therefore, evaluating the both sides at $1$ we get
\[
p^{\#}_{\rho} (\lambda) =\bigl(D_\rho
\chi^\lambda\bigr) (1).
\]
Next, taking the expectation of the both sides with respect to $M_N$,
we get
\[
\mathbf E_N \bigl(p_{\rho}^{\#}\bigr)=
(D_\rho\chi_N ) (1).
\]

Finally, under the specialization of the coefficients of the operator
$D_\rho$
at the point $1\in\operatorname{Mat}(N,\mathbb C)$ this operator
simplifies and
turns into the operator in \eqref{difOp1}.
\end{pf*}

\subsection{Geometric interpretation of signatures}
\label{graph-sign}

Let us depict signatures $\lambda\in\mathbb{GT}_N$ in the way shown
in Figure~\ref{signature}. This figure explains how to assign to $\lambda$ a continuous
piecewise linear function $w_{\lambda}(x)$ (bold line in the figure).

Formally, $w_{\lambda}(x)$ is uniquely determined by the following properties:
\begin{itemize}
\item$w'_{\lambda}(x)$ may have jump discontinuities only at points
$n\in\mathbb Z$ of the
$x$-axis;

\begin{figure}[b]

\includegraphics{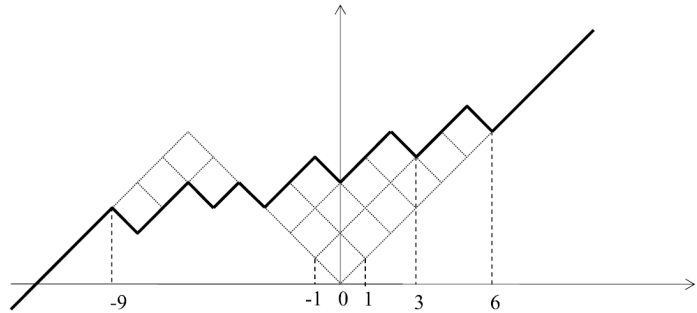}

\caption{A piecewise linear function corresponding to the signature
$\lambda= (6,4,2,0,-1,-3)$.}
\label{signature}
\end{figure}

\item$w'_{\lambda}(x)=\pm1$ for $x\notin\mathbb Z$;

\item$w_{\lambda} (x)=x$ for $x\ge\lambda_1$ and $w_{\lambda}
(x)=x+2N$ for
$x\le\lambda_N-N$, so that $w'_\lambda(x)=1$ outside $[\lambda
_N-N,\lambda_1]$;

\item inside $(\lambda_N-N,\lambda_1)$, there are exactly $N$ unit intervals
$(n,n+1)$ where $w'_\lambda(x)=-1$, and these are those with the midpoints
$\lambda_i-i+\frac{1}2$, $i=1,\ldots,N$.
\end{itemize}

In particular, the function $w_0(x)$ corresponding to the signature
$(0,\ldots,0)\in\mathbb{GT}_N$ has exactly two derivative jumps, at the
points $x=-N$ and
$x=0$.

We regard $w_\lambda$ as the \emph{shape} of $\lambda$. Note that the
part of
the graph of $w_\lambda$ above (resp., below) the broken line $w_0$
visualizes the diagram $\lambda^+$ (resp., $\lambda^-$); see the
\hyperref[sec1]{Introduction} for the definition of $\lambda^\pm$.

We also need the function $\frac{1}Nw_\lambda(Nx)$, which describes the
\emph{scaled shape} of $\lambda$.

Next, recall  definition \eqref{eq6} of the probability measure on
$\mathbb
R$ associated with~$\lambda$:
\[
\mu_{\lambda} := \frac{1}{N} \sum_{i=1}^N
\delta \biggl( \frac{\lambda_i
- i +
{1}/{2}}{N} \biggr),
\]
where $\delta(x)$ denotes the Dirac measure at $x$. Clearly, $\lambda$ is
uniquely determined by~$\mu_\lambda$.

We are going to show that the concentration of random measures $\mu
_\lambda$
implies the concentration of the scaled shapes.

\begin{proposition}\label{prop2}
Assume that for every $N=1,2,\ldots$ we are given an ensemble of random
signatures $\lambda=\lambda(N)$ distributed according to a probability measure
on $\mathbb{GT}_N$. Next, let us assume that, as $N\to\infty$, the
corresponding random measures $\mu_{\lambda}$ weakly converge, in probability,
to a nonrandom probability measure $\sigma$ with support in a bounded interval
$[a,b]\subset\mathbb R$.
\begin{longlist}[(ii)]
\item[{(i)}] The limit measure $\sigma$ is absolutely continuous with
respect to
Lebesgue measure on $\mathbb R$ and so has a density $p(x)$ vanishing outside
$[a,b]$.

\item[{(ii)}] The random functions $\frac{1}N w_\lambda(Nx)$ uniformly
converge in
probability to a nonrandom function $w(x)$, uniquely determined by the
following three properties: $w(x)=x$ for $x>b$, $w(x)=x+2$ for
$x<a$ and
$w'(x)=1-2p(x)$ almost everywhere on $[a,b]$.
\end{longlist}
\end{proposition}

\begin{pf}
(i) The assumption of the proposition means that for any bounded continuous
function $f(x)$ on $\mathbb R$,
%
\begin{equation}
\label{WeakConvergence} \lim_{N \to\infty} \langle f, \mu_{\lambda
}
\rangle= \langle f, \sigma\rangle\qquad \mbox{in probability},
\end{equation}
where the angular brackets denote the pairing between functions and measures.

Let us assume that $f$ is compactly supported and nonnegative. By the very
definition of $\mu_\lambda$,
\[
\langle f, \mu_{\lambda}\rangle=\frac1N\sum_{i=1}^Nf\biggl(
\lambda_i-i+\frac12\biggr)\le \frac1N\sum_{n\in\mathbb
Z}f\biggl(n+
\frac12\biggr).
\]
Since the last expression is the Riemann sum for the integral of $f$ against
Lebesgue measure, passing to a limit as $N\to\infty$, we see that
$\langle f,
\sigma\rangle$ is bounded from above by that integral. If follows that
$\sigma$
has a density $p(x)$ with respect to Lebesgue measure and, moreover,
$p(x)\le1$
almost everywhere.

(ii) Let us define an auxiliary piecewise linear function $\tilde
w_{\lambda}(x)$ by
\[
\tilde w_{\lambda} (x) := x+2 \bigl( 1- \mu_{\lambda} \bigl((-\infty;x]\bigr)
\bigr) = x+ 2 \mu_{\lambda} \bigl((x;+\infty)\bigr).
\]
It readily follows that for $x$ such that $w'_{\lambda} (Nx)=1$ we have
$\frac{1}{N} w_{\lambda} (Nx) = \tilde w_{\lambda}(x)$, and for all $x$,
%
\begin{equation}
\label{estimateAux}\biggl |\tilde w_{\lambda} (x) - \frac{1}{N}
w_{\lambda} (Nx)\biggr| \le\frac{1}{N}.
\end{equation}

Let us define $w(x)$ as the primitive function of $1-2p(x)$ such that
$w(x) =
x$ for $x \gg0$. By virtue of \eqref{WeakConvergence} and claim (i),
$\mu_{\lambda} (\mathbb R \setminus[a,b])$ converges in probability to
0 as
$N\to\infty$. The uniform convergence of $\tilde w_{\lambda}(x)$ to $w(x)$
outside of $[a,b]$ directly follows from this fact. Equation
\eqref{estimateAux} implies that the functions $\frac{1}{N} w_{\lambda} (Nx)$
also uniformly converge to $w(x)$ outside of $[a,b]$.

The definition of $\tilde w_{\lambda} (x)$ and the convergence of
$\mu_{\lambda} (\mathbb R \setminus[a,b])$ to 0 implies that for any bounded
continuous function $f$ we have
\[
\lim_{N \to\infty} \int_a^b f(x)
\tilde w_{\lambda} (x) \,dx = \int_a^b f(x)
w(x) \,dx \qquad\mbox{in probability}.
\]
Using \eqref{estimateAux}, we obtain
%
\begin{equation}
\label{convShapes} \lim_{N \to\infty} \int_a^b
f(x) \frac{1}{N} w_{\lambda} (Nx) \,dx = \int_a^b
f(x) w(x) \,dx \qquad\mbox{in probability}.
\end{equation}
Note that the functions $\frac{1}{N} w_{\lambda} (Nx)$ and $w(x)$ are Lipschitz
functions with Lipschitz constant 1, and for such functions, the
convergence of
the integrals \eqref{convShapes} with an arbitrary continuous test
function $f$
on $[a,b]$ is equivalent to the uniform convergence on $[a,b]$ (see, e.g.,
\cite{IvaOls}, Lemma~5.7). This completes the proof.
\end{pf}

\subsection{Convergence of random measures}
\label{convMe}
In this subsection, we prove a technical lemma about convergence of
random measures.

Let $\{ X_{i,j} \}_{i=1,2, \ldots; j=1,2, \ldots,i}$ be a set of random
variables. Let
\[
\nu_N := \frac{1}{N} \sum_{i=1}^N
\delta( X_{i,N} )
\]
be a (random) measure on $\mathbb R$. Assume that the following
conditions hold:
%
\begin{eqnarray}
\label{meCond1} \lim_{N \to\infty} \mathbf E \int x^k
\nu_N (d x) &=& \mathbf a_k,\qquad k=1,2,3, \ldots,
\\
\label{meCond2} \lim_{N \to\infty} \mathbf E \biggl( \int
x^k \nu_N (d x) \biggr)^2& =& \mathbf
a_k^2,\qquad k=1,2,3, \ldots.
\end{eqnarray}
Also assume that there exists a constant $C>0$ such that
%
\begin{equation}
\label{meCond3} \mathbf a_k < C^k,\qquad k=1,2,3, \ldots.
\end{equation}

\begin{lemma}
\label{lemma-shod}
Let $\{ X_{i,j} \}_{i=1,2, \ldots; j=1,2, \ldots,i}$ be a set of random
variables such that conditions \eqref{meCond1}--\eqref{meCond3} hold.
Then there exists a measure $\nu$ such that
\[
\int x^k \nu(d x) = \mathbf a_k, \qquad k=1,2,3, \ldots,
\]
and we have
\[
\lim_{N \to\infty} \nu_N = \nu\qquad \mbox{weakly; in
probability}.
\]

In greater detail, for any bounded continuous $f$ we have
\[
\lim_{N \to\infty} \int f\, d \nu_N = \int f\, d \nu,\qquad
\mbox{in probability}.
\]
\end{lemma}

\begin{pf}
We follow \cite{AGZ}, Section~2.1.2.

Define a (deterministic) measure $\bar\nu_N$ on $\mathbb R$ by its
values on test functions via
\[
\int f d \bar\nu_N := \mathbf E \int f \,d \nu_N\qquad
\mbox{for any bounded continuous $f$}.
\]

It follows from the Chebyshev inequality that for any $B >1$ we have
\[
\mathbf P \biggl( \int x^k \mathbf1_{|x|>B}
\nu_N (d x) > \varepsilon \biggr) \le\frac{1}{\varepsilon} \mathbf E \int
x^k \mathbf1_{|x|>B} \nu_N (d x) \le
\frac{\mathbf E \int x^{2k} \nu_N (d x)}{\varepsilon B^k}.
\]

Conditions \eqref{meCond1} and \eqref{meCond3} imply
\[
\limsup_{N \to\infty} \mathbf P \biggl( \int x^k
\mathbf1_{|x|>B} \nu _N (d x) > \varepsilon \biggr) \le
\frac{\mathbf a_{2k}}{\varepsilon B^k} \le \frac{(C^2)^k}{\varepsilon B^k}.
\]

Note that for any $K > k$ we have
\begin{eqnarray*}
\limsup_{N \to\infty} \mathbf P \biggl( \int x^k
\mathbf1_{|x|>B} \nu _N (d x) > \varepsilon \biggr) &\le&\limsup
_{N \to\infty} \mathbf P \biggl( \int x^{2K}
\mathbf1_{|x|>B} \nu_N (d x) > \varepsilon \biggr)\\
& \le&
\frac
{C^{4K}}{\varepsilon B^{2K}}.
\end{eqnarray*}

Choosing $B=C^2+1$ and letting $K$ to infinity, we have
%
\begin{equation}
\label{eq0} \limsup_{N \to\infty} \mathbf P \biggl( \int
x^k \mathbf1_{|x|>B} \nu _N (d x) > \varepsilon
\biggr) = 0.
\end{equation}
Therefore, we obtain
\[
\lim_{N \to\infty} \mathbf E \int x^k
\mathbf1_{[-B;B]} \nu_N (d x) = \mathbf a_k,\qquad
k=1,2,3, \ldots.
\]

Since the unit ball in $(C[-B;B])^*$ is weakly compact, the sequence
$\bar\nu_N$ converges (weakly) to a probability measure $\nu$ with
support in $[-B;B]$, and we have
\[
\int x^k \nu(d x) = \mathbf a_k.
\]

Note that \eqref{meCond1} and \eqref{meCond2} imply that the sequence
$\int x^k \,d \nu_N$ converges to $\mathbf a_k$ in probability.

Let $f(x)$ be a continuous bounded function. The Weierstrass theorem
implies that for any $\delta>0$ there exists a polynomial $Q_{\delta
}(x)$ such that
\[
\sup_{x \in[-B;B]} \bigl|Q_{\delta} (x) - f(x)\bigr| < \delta/10.
\]
Then
\begin{eqnarray*}
&&\mathbf P \biggl( \biggl\llvert \int f(x) \nu_N (d x) - \int f(x)
\nu(d x) \biggr\rrvert > \delta \biggr) \\
&&\qquad\le\mathbf P \biggl( \biggl\llvert \int
Q_{\delta} (x) \nu_N(d x)- \int Q_{\delta} (x) \nu(d x)
\biggr\rrvert > \delta/4 \biggr)
\\
&&\qquad\quad{}+ \mathbf P \biggl( \biggl\llvert \int Q_{\delta} (x)
\mathbf1_{|x|>B} \nu _N( d x) \biggr\rrvert > \delta/4
\biggr).
\end{eqnarray*}
The first term converges to zero due to the convergence in probability
of $\int x^k \nu_N( d x)$ to $\mathbf a_k$, and the second term
converges to zero due to \eqref{eq0}. This completes the proof of the lemma.
\end{pf}

\section{Main result and discussion}\label{sec3}

\subsection{The main result}
\label{statement}
In this section, we state the main result of this paper.\vadjust{\goodbreak}

Recall that elements $\omega\in\Omega$ parameterize extreme
characters of the group $U(\infty)$. We consider a sequence $\omega(N)
\in\Omega$ depending on (growing) integer $N$.
Let $\chi^{\omega(N)}$ be the extreme character of $U(\infty)$
corresponding to $\omega(N)$, and let $M_N$ be the probability measure
on $\mathbb{GT}_N$ determined by this character (see Section~\ref{gelfand-tsetlin}).
Let $\lambda^{(N)} \in\mathbb{GT}_N$ be a random
signature distributed according to $M_N$, and let
\[
\mu^{(N)} := \mu_{\lambda^{(N)}} = \frac{1}{N} \sum
_{i=1}^N \delta \biggl( \frac{\lambda_i^{(N)} - i+1/2}{N} \biggr)
\]
be the random measure on $\mathbb R$ associated with $\lambda^{(N)}$ (see
Section~\ref{graph-sign}). We are interested in the limit behavior of this
random measure.

Let $\Phi^{\omega(N)} (z)$ be the Voiculescu function depending on parameters
$\gamma^{\pm} (N)$, $\{\alpha_i^{\pm} (N) \}$, $\{ \beta_j^{\pm} (N) \}$,
see \eqref{eq4}. Consider the following condition on these
sequences of parameters.

\textit{Main condition}. Assume that for some $\varepsilon>0$ the analytic
function\break $\log\Phi^{\omega(N)}(z+1)$ uniformly converges to an analytic
function $P(z)$ on $\{z \in\mathbb C  | |z| \le\varepsilon\}$:
%
\begin{equation}
\label{cond-unit} \lim_{N \to\infty} \frac{1}{N} \bigl(\log
\Phi^{\omega(N)}(z+1) \bigr) = P(z).
\end{equation}

By $t_i$, $i\in\mathbb N$, we denote the coefficients of the Taylor
series for $P'(z)$:
%
\begin{equation}
\label{deriv} P'(z) =: t_1 + t_2 z +
t_3 z^2 + \cdots.
\end{equation}

It is convenient for us to also formulate a stronger condition that describes
more explicitly how the parameters $\gamma^{\pm} (N)$, $\{\alpha_i^{\pm
} (N)
\}$, $\{ \beta_j^{\pm} (N) \}$ can change.

\textit{Sufficient condition.}
Let
\[
\mathcal A_N^{\pm} := \frac{1}{N} \sum
_{i=1}^{\infty} \delta\bigl(\alpha _i^{\pm}
(N) \bigr), \qquad\mathcal B_N^{\pm} := \frac{1}{N} \sum
_{i=1}^{\infty} \delta \bigl(
\beta_i^{\pm} (N) \bigr),
\]
be measures on $\mathbb R$.

We say that a sequence $\omega(N)$ satisfies the sufficient condition
if there exist limits
%
\begin{equation}
\label{sufCond1} \lim_{N \to\infty} \frac{\gamma^{\pm} (N) }{N}
\end{equation}
and
%
\begin{equation}
\label{sufCond2} \lim_{N \to\infty} \mathcal A_N^{\pm}
= \mathcal A^{\pm},\qquad \lim_{N \to\infty} \mathcal
B_N^{\pm} = \mathcal B^{\pm}\qquad \mbox{weak
convergence},
\end{equation}
for some finite measures $\mathcal A^{\pm}$, $\mathcal B^{\pm}$ on
$\mathbb{R}$
with compact support. Moreover, we require that there exist positive constants
$C_1$, $C_2$ such that
%
\begin{equation}
\label{sufCond3} \bigl|\alpha_i^{\pm}(N)\bigr| < C_1,\qquad\bigl |
\beta_i^{\pm} (N)\bigr| < C_1\qquad \mbox{for all $i
\ge1$},
\end{equation}
and the number of nonzero parameters $\alpha_i^{\pm}(N)$ and $\beta
_i^{\pm}
(N)$ is less than $C_2 N$.

For example, let $\alpha_1 = \cdots= \alpha_N = \alpha$, where
$\alpha>0$ is a fixed constant, and all other Voiculescu's parameters
are equal
to 0. It is clear that this sequence of parameters satisfies the sufficient
condition \eqref{sufCond1}--\eqref{sufCond3}. Another example is given by
$\gamma^+=\gamma N$, where $\gamma>0$ is a fixed constant, and all other
Voiculescu's parameters are equal to 0. More examples can be found in the
\hyperref[app]{Appendix}.

\begin{proposition}
Let $\{ \omega(N) \}_{N \ge1}$ be a sequence of points in $\Omega$.
Assume it satisfies the sufficient condition \eqref{sufCond1}--\eqref
{sufCond3}. Then it also satisfies the main condition \eqref{cond-unit}.
\end{proposition}

\begin{pf}
We have (omitting the dependence on $N$ in notation)
%
\begin{eqnarray}
\label{log} &&\frac{1}{N} \log\Phi^{\omega(N)} (z+1) \nonumber\\
&&\qquad=
\frac{\gamma^+}{N} + \frac
{\gamma^-}{N} \biggl( \frac{1}{z+1} -1 \biggr)
\nonumber
\\[-8pt]
\\[-8pt]
\nonumber
&&\qquad\quad{}+\frac{1}{N} \biggl( \sum_{i\ge1} \log\bigl(1+
\beta_i^+ z\bigr) - \sum_{i\ge1} \log
\bigl(1- \alpha_i^+ z\bigr)\\
&&\hspace*{30pt}\qquad\quad{} + \sum_{i \ge1}
\log \biggl( 1 -\frac{\beta_i^- z}{1+z} \biggr) - \sum_{i
\ge1}
\log \biggl( 1+ \frac{\alpha_i^- z}{1+z} \biggr) \biggr).
\nonumber
\end{eqnarray}
Conditions \eqref{sufCond2} and \eqref{sufCond3} imply that there exist limits
\[
\frac{1}{N} \sum_i \mathbf
s_i^k \qquad\mbox{for all $k \ge0$},
\]
where $\mathbf s_i$ is equal to $\alpha_i^{\pm}$ or $\beta_i^{\pm}$.

This fact and condition \eqref{sufCond1} imply that the Taylor
coefficients of \eqref{log} converge to some limiting coefficients
$t_i$, and the power series determined by these $t_i$ converges in a
neighborhood of 0 (because the supports of $\mathcal A^{\pm}$ and
$\mathcal B^{\pm}$ are compact). Condition \eqref{sufCond3} implies
that this convergence is uniform.
\end{pf}

From now on we assume that $\omega= \omega(N)$ satisfies the main
condition \eqref{cond-unit}.

Let
\[
Q(z) = 1 + z(1+z) \bigl(t_1 + t_2 z + t_3
z^2 + \cdots\bigr)
\]
be a formal power series depending on coefficients $t_1, t_2, \ldots.$
Define a formal power series $v_0(z)$ via
\[
v_0(z) := \biggl( \frac{z}{Q(z)} \biggr) ^{(-1)},
\]
where in the right-hand side the formal inversion of power series is
used.\footnote{Here and below we consider formal power series of the form
\[
z+ a_2 z^2 + a_3 z^3 + \cdots,\qquad
a_i \in\mathbb R.
\]
It is well known that such a series has a unique inverse (with respect
to composition) of this form. For example, if $A(z) = \sum_{i=1}^{\infty
} z^i$ then $A^{(-1)} (z) = \sum_{i=1}^{\infty} (-1)^{i-1} z^i$.}

Let
%
\begin{equation}
\label{Sz} S(z) := \log\bigl(1+ v_0 (z) \bigr) = z+
m_1 z^2 + m_2 z^3 + \cdots.
\end{equation}
Later on (see Section~\ref{planProof}) we will prove that there exists
a unique probability measure on $\mathbb R$ with moments $\{1,m_1, m_2,
\ldots\}$. Denote this measure by $\sigma$.

The main result of this paper is
the following.

\begin{theorem}
\label{mainTh} Let $\{ \omega(N) \}_{N \ge1}$ be a sequence of points
in $\Omega$
satisfying condition \eqref{cond-unit}. Then
\[
\lim_{N \to\infty} \mu^{(N)} = \mathcal\sigma \qquad\mbox{weak
convergence in probability.}
\]
Equivalently, for any bounded continuous function $f$ we have
\[
\lim_{N \to\infty} \int f d \mu^{(N)} = \int f d \sigma\qquad
\mbox{in probability.}
\]
\end{theorem}

The proof of this theorem is given in Section~\ref{proofTh}.

The density of $\sigma$ (which is well defined by virtue of Proposition~\ref{prop2}) can be obtained from the known Stieltjes transform
$S(1/z)$ [see
\eqref{Sz}] with the use of standard methods of complex analysis; see
also the end of the next subsection.

\subsection{A heuristic derivation of the limit shape}
\label{DetPr}

In this section, we sketch an argument which shows how one can compute the
measure $\sigma$ (see Theorem~\ref{mainTh}) via determinantal point processes.
This yields the correct formula but not a complete proof, because the very
existence of the concentration remains unclear. Our proof of Theorem~\ref{mainTh} is obtained in a very different way (see Section~\ref{proofTh}).

In \cite{BorKuan}, it was shown that the correlation functions of the
random point configuration $(\lambda_1 -1, \lambda_2-2, \ldots)$
corresponding to the restriction of the extreme character of $U(\infty
)$ with Voiculescu function $\Phi^{\omega} (z)$ to $U(N)$ have
determinantal structure (necessary definitions can be found, e.g., in
\cite{BorKuan}). The correlation kernel of this process has the
following form:
\[
K(x,y) = \frac{-1}{4 \pi^2} \oint\oint\frac{\Phi^{\omega
}(u^{-1})}{\Phi^{\omega}(w^{-1})} \frac{u^x (1-u)^N}{w^{1+y} (1-w)^N}
\frac{du \,dw}{u-w},
\]
where the $u$-contour is a counterclockwise oriented circle with center
0 and radius $\varepsilon\ll1$, and the $w$-contour is a
counterclockwise oriented circle with center 1 and radius $\delta\ll1$.

If one already knows that the random point process $(\lambda_1 -1,
\lambda_2-2, \ldots)$ satisfies the law of large numbers type theorem,
then it is natural to assume that the density of the limit measure is
equal to the limit of the diagonal values of the kernel $N^{-1}
K(xN,xN)$ (this is the so-called density function) as $N \to\infty$.

Let us find (informally) the limit of $N^{-1} K(xN,xN)$ as $N \to\infty
$. A useful general approach to asymptotic analysis of such integrals
is the steepest decent method. In order to apply this method, we write
the integrand in the form
\[
\frac{\exp ( N  ( \mathcal S(z) - \mathcal S(w)  )
) }{z-w},
\]
where
\[
\mathcal S(u) := \frac{\log f(u^{-1}) + x \log u + N \log(1-u)}{N}.
\]

Following the logic of \cite{Oko} (see also \cite{BorGor}), we need to
deform the contours of integration in such a way that they pass through
the critical points of $\mathcal S(z)$ which are the roots of
%
\begin{equation}
\label{eq-shape} \frac{1}{N} \biggl( \bigl(\log\Phi^{\omega}
\bigl(z^{-1}\bigr)\bigr)' + \frac{x}{z} -
\frac
{N}{1-z} \biggr) = 0.
\end{equation}
We are interested in the root $z_+=z_+(x)$ which has the positive
imaginary part.

Then the steepest decent method gives the following asymptotics for the
one-dimensional correlation function (cf. \cite{Oko,BorGor}):
%
\begin{equation}
\label{limitDensity} \frac{1}{N} K(xN,xN) \approx\frac{1}{\pi} \arg(z_+),\qquad
N \to\infty.
\end{equation}

Let us apply a change of variable $z = 1/w$; with the use of \eqref
{cond-unit}, equation \eqref{eq-shape} can be written in the form
%
\begin{equation}
\label{eq-shape2} P'(w-1) - \frac{x+1}{w} + \frac{1}{w-1} =
0.
\end{equation}
Let $w_0 = w_0 (x)$ be the complex root of \eqref{eq-shape2} in the
complex upper half-plane.

Recall that the Stieltjes transform of a probability measure
$\hat\mu$ with compact support is given by
\[
\mathrm{Stil}_{\hat\mu} (z) := \int_{\mathbb R} \frac{\hat\mu(dt)}{z - t},\qquad z
\in\mathbb C \setminus\operatorname{supp}(\hat\mu).
\]
Observe that if one denotes the moments of $\hat\mu$ by $1,m_1,m_2,\ldots,$ then $\mathrm{Stil}_{\hat\mu}(z)$ is obtained from the right-hand side of
\eqref{Sz} by the change of variable $z\mapsto z^{-1}$.

The Stieltjes transform can be inverted. For a measure $\hat\mu$ with density
$\hat p (x)$ with respect to the Lebesgue measure, we have
\[
\hat p(x) = \lim_{\varepsilon\to0} \frac{1}{\pi} \Im
\bigl(\mathrm{Stil}_{\hat\mu} (x + i \varepsilon)\bigr).
\]
Assume that $w_0=w_0(x)$ is a real-analytic function. Using the analytic
continuation, one can view $w_0=w_0(x)$ as a complex analytic function.
It is
natural to think that for the principal branch of the function $\log
(x)$ we
have
\[
\lim_{\varepsilon\to0} \Im\log\bigl(w_0(x+ i \varepsilon)
\bigr) = \arg w_0(x),\qquad x\in\mathbb R.
\]
Note also that $\arg(z_+(x)) = \arg(w_0(x))$. Therefore, it is natural to
assume that the Stieltjes transform of the limit measure is equal to
$\log
(w_0(x))$.

Let $v_0(z)$ be the formal power series defined in Section~\ref{statement}. Note that the series $y_0 = v_0(1/z)$ solves the following
equation:
%
\begin{equation}
\label{eq-shape3} z = \frac{1}{y} + (1+y) P'(y).
\end{equation}

Equations \eqref{eq-shape2} and \eqref{eq-shape3} imply that the formal
power series $v_0(1/z)$ satisfies the same equation as $w_0-1$. Thus,
the result stated in Theorem~\ref{mainTh} coincides with the heuristic
answer coming from the determinantal processes.

\subsection{Markov dynamics on two-dimensional arrays}\label{sec3.3}

This subsection details the relation of the present work to random
growth of surfaces in $(2+1)$-dimensions. This connection served as our
original motivation, but it is not necessary for understanding the rest
of the paper, and thus the reader should feel free to omit it.

Consider a two-dimensional triangular array of particles
\[
\mathcal W = \bigl\{ \bigl\{ x_k^m \bigr
\}_{m=1, \ldots, \infty; k = 1, \ldots, m} \subset\mathbb Z^{{n(n+1)}/{2}} | x_k^{m+1}
\ge x_k^m > x_{k+1}^{m+1} \bigr\};
\]
we interpret the number $x_k^m$ as the position of the particle with
label $(k,m)$.

For any $N$, the extreme character $\omega(N)$ is determined by a set of
Voicu\-lescu's parameters $\alpha^{\pm} (N)$, $\beta^{\pm} (N)$, and
$\gamma^{\pm} (N)$ (below we omit the dependence on $N$ in notation). Suppose
that the number of parameters of types $\alpha^{\pm}$, $\beta^{\pm}$ is finite
and equals $T=T(N)$. Let us enumerate these parameters by the
numbers $1, \ldots, T$ in an arbitrary way. We interpret this enumeration
dynamically as follows: At time 1, we take only the first parameter; at
time 2,
the second parameter is added, etc. Let $\chi_a$ be the extreme
character of
$U(\infty)$ determined by the first $a$ parameters in our ordering, $1
\le a
\le T$. This character gives rise to a probability measure on $\mathbb{GT}_N$;
denote it by $\mu^{(N)}_a$.

It turns out that the measure $\mu^{(N)}_a$ can be obtained in the following
way. One can define (see \cite{BorFer}, Section~2.6) a discrete time Markov
dynamics on the triangular arrays as above with the following property:
For any
$N \ge1$, at each time $a$ the distribution of the vector $\{ x^N_k - N
\}_{k=1, \ldots, N}$ coincides with the distribution of $\{ \lambda_i - i
\}_{i\ge1}$, where $\{ \lambda_i \}$ are the coordinates of the random
signature distributed according to the measure $\mu^{(N)}_a$. In particular,
for $a=T$ the distribution of the $N$th level of the array coincides
with the
measure $\mu^{(N)}$. Parameters $\gamma^{\pm}$ can also be realized
under a
similar Markov dynamics with continuous time (see \cite{BorFer}, Section~1).

An important feature of these Markov processes is the locality of
interactions between the particles---the behavior of each individual
particle is only influenced by particles whose
coordinates differ at most by 1 from those of the chosen one.

The evolution of the whole array of particles can be fully encoded by
the \textit{height function} $h\dvtx  \mathbb R \times\mathbb R_{\ge1}
\times\{1,2, \ldots, T \} \to\mathbb Z_{\ge0}$ defined by
\[
h(\mathbf x, \mathbf y, a) = \# \bigl\{ k | x_k^{[\mathbf y]} (a) >
\mathbf x \bigr\},
\]
where $x_k^m (a)$ stands for the position of the particle $x_k^m$ at
the time $a$.

Suppose now that a sequence of characters $\omega(N)$ satisfies the
general condition~\eqref{cond-unit} with a function $P_0(z)$. Let us
fix a large $N$; at this stage, we have a certain set of parameters $\{
\alpha^{\pm}, \beta^{\pm} \}$. We want to add to these parameters
another set of parameters satisfying condition~\eqref{cond-unit}. For
simplicity, we consider six special cases: adding $t N$ parameters of
one of the possible types $\alpha^{\pm}$, $\beta^{\pm}$, or increasing
$\gamma^{\pm}$ by $t N$. Then the function $P(z)$ describing such a
model can be written in the form
\[
P(z) = P_0(z) + t F(z),
\]
where $F(z)$ is determined by the choice of one of the special cases
mentioned above. Let
\[
h(x,y,T + t) = \lim_{N \to\infty} \frac{\mathbf E h ([xN], [yN], (T+ t)N)}{N}
\]
be the limiting height function.

The plot of the height function with a fixed 3rd coordinate can be
viewed as
a random two-dimensional surface in $\mathbb R^3$.
As was mentioned above, the growth of the height function can be
realized as a result of Markov dynamics with local interactions.
The theory of hydrodynamic limits of random growth models allows one to
predict the type of the modification of the limit shape when we add
parameters with the use of local Markov dynamics. Namely, one can
expect that the limit height function obeys an evolution equation of
the form
\[
\frac{\partial h(x,y,T + t)}{\partial t} = \mathcal F \biggl( \frac
{\partial h(x,y,T + t)}{\partial x}, \frac{\partial h(x,y,T +
t)}{\partial y}
\biggr),
\]
where $\mathcal F$ is a function of two variables uniquely determined
by the function $F$ (or, equivalently, by the type of parameters that
we add).

Let us verify that the limit measure coming from Theorem~\ref{mainTh}
satisfies such an equation.
We use the answer in the form given in Section~\ref{DetPr}, namely, let
\begin{eqnarray*}
S(z) &=& P \biggl(\frac{1}{z}-1 \biggr) + x \log z + y \log(1-z) \\
&=&
P_0 \biggl( \frac{1}{z}-1 \biggr) + t F \biggl(
\frac{1}{z}-1 \biggr) + x \log z + y \log(1-z).
\end{eqnarray*}
Then the density of the limit measure is equal to $\frac{1}{\pi} \arg
(z_+ (x,y,t))$, where $z_+$ is the root of the equation $S'(z)=0$ lying
in the upper half-plane. Hence,
\[
h(x,y,T + t) = - \frac{1}{\pi} \Im\bigl(S\bigl(z_+ (x,y,t)\bigr)\bigr).
\]

Differentiating this equality and taking into account that $S'(z_+)=0$,
we obtain
\begin{eqnarray*}
\frac{\partial h(x,y,T + t)}{\partial x} &=& -\frac{1}{\pi} \Im\bigl( \log (z_+)\bigr) = -
\frac{1}{\pi} \arg(z_+),
\\
\frac{\partial h(x,y,T + t)}{\partial y} &=& -\frac{1}{\pi} \Im\bigl( \log(1 - z_+)\bigr) = -
\frac{1}{\pi} \arg(1 - z_+),
\\
\frac{\partial h(x,y,T + t)}{\partial t} &=& - \frac{1}{\pi} \Im \biggl( F \biggl(
\frac{1}{z_+}-1 \biggr) \biggr).
\end{eqnarray*}
Note that the arguments of $z_+$ and $1-z_+$ uniquely determine the
complex number $z_+$ with a positive imaginary part. Therefore, the
function $F$ and the derivatives of $h(x,y,T +t)$ with respect to $x$
and $y$ uniquely determine the derivative of $h(x,y,T +t)$ with respect
to $t$.

\begin{figure}

\includegraphics{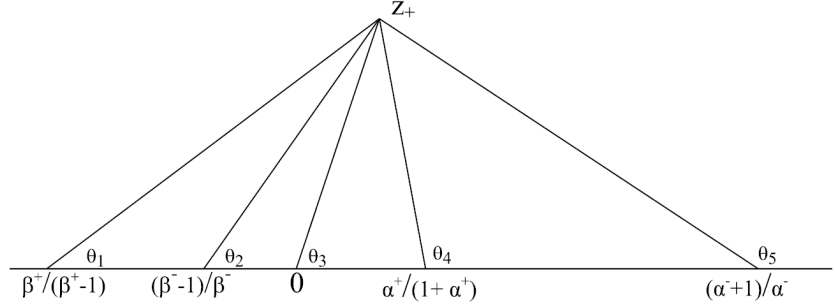}

\caption{Angles which determine the growth of the limit shape.}\label{Ugly}
\end{figure}

When we add equal $\alpha^+$ parameters we have $F(z) = - \log(1 -
\alpha^+ z)$. After computations, we obtain
\[
\frac{\partial h(x,y,T + t)}{\partial t} = \frac{1}{\pi} \biggl( \arg \biggl( z_+ -
\frac{\alpha^+}{1+\alpha^+} \biggr) - \arg(z_+) \biggr) = \frac{\theta_4 - \theta_3}{\pi},
\]
where the angles $\theta_i$ are shown in Figure~\ref{Ugly}.
%

Analogous computations for five other cases (equal $\beta^+$'s, $\alpha
^-$'s, $\beta^-$'s and the growth of $\gamma^+$ or $\gamma^-$) show
that (similar computations were performed in \cite{Chen}), respectively,
\begin{eqnarray*}
\frac{\partial h(x,y,T + t)}{\partial t} &=& \frac{1}{\pi} \biggl(- \arg \biggl(z_+ +
\frac{\beta^+}{1-\beta^+} \biggr) - \arg(z_+) \biggr) = \frac{- \theta_1 - \theta_3}{\pi},
\\
\frac{\partial h(x,y,T + t)}{\partial t}& =& \frac{1}{\pi} \biggl( \arg \biggl( z_+ -
\frac{1+\alpha^-}{\alpha^-} \biggr) - \pi \biggr) = \frac
{\theta_5 - \pi}{\pi},
\\
\frac{\partial h(x,y,T + t)}{\partial t} &=& - \frac{1}{\pi} \biggl( \arg \biggl( z_+ +
\frac{1- \beta^-}{\beta^-} \biggr) \biggr) = \frac{\theta
_2}{\pi},
\\
\frac{\partial h(x,y,T + t)}{\partial t} &=& - \frac{\gamma^+}{\pi} \Im \biggl( \frac{1}{z_+}
\biggr),
\\
\frac{\partial h(x,y,T + t)}{\partial t} &=& - \frac{\gamma^-}{\pi} \Im(z_+).
\end{eqnarray*}

\subsection{A convolution of measures}\label{sec3.4}

Let $\mathcal M$ be the space of probability measures that can be
obtained as limit measures $\sigma$ from Theorem~\ref{mainTh}. For $\nu
\in\mathcal M$, let $P'_{\nu}$ be the function defined in \eqref
{cond-unit}, and let $S_{\nu}$ be the generating function of moments
defined in \eqref{Sz}.
These two functions uniquely determine each other by
%
\begin{equation}
\label{linVmera} S_{\nu}(z) = \log \biggl( 1 + \biggl( \frac{z}{1 + z(1+z) P'_{\nu}(z)}
\biggr)^{(-1)} \biggr)
\end{equation}
and
%
\begin{equation}
\label{meraVlin} P'_{\nu} (z) = \frac{1}{(1+z) (\exp(S_{\nu} (z))-1)^{(-1)}} -
\frac{1}{z(1+z)}.
\end{equation}

Let $\chi^1$ and $\chi^2$ be two extreme characters of $U(\infty)$. Consider
the product of these characters, which is also an extreme character of
$U(\infty)$:
\[
\chi^{1,2} (U) := \chi^1 (U) \chi^2 (U), \qquad U \in
U(\infty).
\]
It is natural to think that this operation corresponds to a tensor
product of
representations of $U(\infty)$ determined by the characters $\chi^1$ and
$\chi^2$ (although these are infinite-dimensional objects and one needs to
explain what that means).

Assume that $\chi^1_N$ and $\chi^2_N$ are sequences of extreme
characters of
$U(\infty)$ satisfying condition \eqref{cond-unit}. By Theorem~\ref{mainTh},
there are limit measures $\sigma_1$ and $\sigma_2$ corresponding to these
sequences. Then the sequence $\chi^1_N \chi^2_N$ also satisfies
\eqref{cond-unit}; let $\sigma_{1,2}$ be the limit measure for this sequence.
Note that
%
\begin{equation}
\label{linear} P'_{\sigma_{1,2}} (z) = P'_{\sigma_1}
(z) + P'_{\sigma_2} (z).
\end{equation}

Thus, these formulas allow to define a natural operation of ``quantized free
convolution'' for measures $\sigma_1, \sigma_2 \in\mathcal M$; the
result of
convolution is $\sigma_{1,2} \in\mathcal M$. The measure $\sigma
_{1,2}$ is
completely determined by equations \eqref{meraVlin}, \eqref{linear} and
\eqref{linVmera}.

Special cases considered in the \hyperref[app]{Appendix} can serve as examples of this
convolution.
In particular, the limit measures for one-sided Plancherel characters with
parameter $\gamma N$ or the characters corresponding to $a N$ parameters
$\alpha_j^+ \equiv1$, form one-parameter subgroups with respect to this
convolution.

This operation of convolution can be defined by the same
formulas for a more general class of measures; the setting for such a
generalization is as follows. Let $T_{\lambda_1}$ and $T_{\lambda_2}$,
$\lambda_1, \lambda_2 \in\mathbb{GT}_N$, be two irreducible
representations of
$U(N)$. Let us consider the Kronecker tensor product $T_{\lambda
_1}\otimes
T_{\lambda_2}$ and decompose it onto irreducible representations. As $N
\to
\infty$, under appropriate scaling regime one can prove a law of large numbers
type theorem for this decomposition; see \cite{BufGor}.

For the first time a similar problem was considered by Biane
\cite{Biane}; the resulting operation on measures was the free convolution.
However, we consider a different scaling, and in our situation the resulting
operation is not the free convolution (see \cite{BufGor}, Section~1, for more
details). In fact, for a certain degeneration turning the branching of
signatures in the Gelfand--Tsetlin graph into the branching of eigenvalues
(describing the eigenvalues of \emph{corners} of Hermitian matrices),
which corresponds to the degeneration of the Gelfand--Tsetlin graph to the
``graph'' of spectra of Hermitian matrices, our convolution turns into
the free
convolution. Let us show how this happens.

Let $R_{\nu} (x)$ be Voiculescu's $R$-function of the measure $\nu$
(see, e.g., \cite{NicSpe}). Then it is easy to see that
\[
P'_{\nu} (z) = \frac{1}{1+z} R_{\nu} \bigl(
\log(1+z)\bigr) + \frac{1}{(1+z) \log
(1+z)} - \frac{1}{z(1+z)}.
\]

For the degeneration to the ``graph'' of spectra, we need to consider
measures with homothetically growing supports and for values of
variables that are close to~1. Let $L$ be a large parameter, and let us
change the variable $z= y/L$. The new $R$-function satisfies
\[
R_{\nu} \bigl(\log(1+z)\bigr) = L R_{\tilde\nu} \biggl( L \log
\biggl(1 + \frac
{y}{L} \biggr) \biggr),
\]
where $\tilde\nu$ is the measure arising after the degeneration.

Thus, we have
\[
\frac{P'_{\nu} (z)}{L} = R_{\tilde\nu} (y) + O \bigl(L^{-1}\bigr)
\mathop{\rightarrow}_{L
\to\infty} R_{\tilde\nu} (y).
\]
Therefore, in this limit the linearizing function $P'_{\nu}$ becomes
the $R$-function of a measure, and the tensor product of
representations gives rise to the free convolution.

\section{Proof of Theorem \texorpdfstring{\protect\ref{mainTh}}{3.2}}
\label{proofTh}

In this section, we prove our main result, Theorem~\ref{mainTh}.
Because the proof is
rather long, let us describe first its main ideas.

To establish the existence of a limit shape, we use the method of moments.
Recall that we interpret signatures $\lambda\in\mathbb{GT}_N$ as
certain measures on $\mathbb{Z}$,
so random signatures become random measures. The moments of the random
measures, as well as products of moments, are thus random functions. We
have to
examine the limit of their expectations as $N\to\infty$.

Our key technical tool is the algebra $\operatorname{Sym}^*$ of shifted
symmetric functions.
As explained in Section~\ref{sec2}, elements of $\operatorname{Sym}^*$ can be
converted, via the maps
$\mathrm{pr}_N$, into functions on signatures. We are dealing with two bases in
$\operatorname{Sym}^*$, $\{\mathbf p_\rho\}$ and $ \{ \mathbf p^\#_\rho
\}$. The products of moments that we need
to control are given by the elements of the first basis, whereas the
expectations are initially expressed in terms of the second basis. This
is the
source of the problem, because the transition coefficients between both basis
have a very complicated structure and hardly can be written down explicitly.

Fortunately, we do not need to know the transition coefficients
exactly, because
for our purpose it suffices to compute their large-$N$ asymptotics, so
that we
may drop many asymptotically negligible terms. This allows us to solve the
problem by reducing it to combinatorial analysis of certain special graphs
(Sections~\ref{proof41} and \ref{proof42}). In the process we recover the noncrossing partitions
which make the connection with free probability (see Section~\ref{sec3.4}) less
surprising.

Note that a similar difficulty of transition between two bases in
$\operatorname{Sym}^*$
arose in Kerov's proof of his central limit theorem for the Plancherel measure
(see Ivanov--Olshanski \cite{IvaOls}). However, in our case the limit
regime is
different, the emerging technical problems are more serious, and the required
combinatorial machinery is substantially more sophisticated.

\subsection{Plan of the proof}\label{planProof}


Let us modify the measure $\mu^{(N)}$ by adding atoms of weight $-\frac{1}{N}$
at locations $-\frac{i}{N}$, $i=1,2, \ldots, N$. Let $\tilde\mu^{(N)}$ denote
the resulting signed measure. Note that its total weight equals 0. As
$N \to
\infty$, the negative part of $\tilde\mu^{(N)}$ converges to the
measure with
density $-1$ on the interval $[-1;0]$.

Recall that the functions from $\mathbb A(N)$ (see Sections~\ref{shifted-symmetric-functions} and \ref{functions-on-signatures}) are
defined on $\mathbb{GT}_N$. By \eqref{functions-pk} the functions $p_k$
are the
moments of the measure $\tilde\mu^{(N)}$. We know the limit of the negative
part of $\tilde\mu^{(N)}$; therefore, the information about the limit
of $\{
p_k \}_{k \ge1}$ is sufficient for describing the limit measure $\sigma
$ (see
Section~\ref{statement}).

Let $\tilde\sigma$ be the sum of the measure $\sigma$ and the negative
Lebesgue measure on $[-1;0]$.
Let $\tilde m_k$ be the moments of $\tilde\sigma$. Define a formal
power series $\tilde S (z)$ by
\[
\tilde S (z) = \tilde m_1 z^2 + \tilde m_2
z^3 + \cdots.
\]
It is easy to see that
%
\begin{equation}
\label{Stilt} S(z) = \tilde S(z) + \log(1+z).
\end{equation}


We recall that by $[u^k]{A(u)}$, where $A(u)$ is a formal power series
of the form $a_1 u + a_2 u^2+ \cdots,$ we denote the coefficient of
$u^k$ in $A(u)$.

Recall that the functions $p_k^{\#}$ are defined on $\mathbb{GT}_N$;
see Section~\ref{functions-on-signatures}.

Let $\mathbf E_N$ denote the expectation with respect to $M_N$.

%
\begin{proposition}
\label{41}
For any $k \ge1$, we have
\[
\lim_{N \to\infty} \frac{\mathbf E_N (p_k^{\#})}{N^{k+1}} = \frac
{1}{k+1}
\bigl[u^k\bigr] \bigl(1+t_1 u + t_2
u^2 + \cdots\bigr)^{k+1} =: c_k.
\]
\end{proposition}

Note that the weight of $p_k^{\#}$ equals $k+1$.

The proof of this proposition is given in Section~\ref{proof41}.

\begin{proposition}
\label{42}
For any partition $\rho= (k_1, k_2, \ldots, k_{l(\rho)})$, we have
\[
\lim_{N \to\infty} \frac{\mathbf E_N (p_{\rho}^{\#} )}{N^{k_1 + k_2 +
\cdots+ k_{l(\rho)}+ l(\rho)}} = c_{k_1}
c_{k_2} \cdots c_{k_{l(\rho)}}.
\]
\end{proposition}

The proof of this proposition is given in Section~\ref{proof42}.

Let us recall that $\{ p_{\rho}^{\#} \}_{\rho\in\mathbb Y}$ is a
linear basis
in $\operatorname{Sym}^*$; therefore, these two propositions give us complete
information about expectations of functions from $\mathbb A(N)$. In particular,
these propositions imply
\[
\lim_{N \to\infty} \mathbf E_N (f) = O\bigl(
N^{\operatorname{wt}(f)}\bigr),
\]
where $\operatorname{wt}(f)$ is the weight filtration.

\begin{proposition}
\label{43}
For any $k \ge1$, we have
\[
\lim_{N \to\infty} \frac{\mathbf E_N (p_k)}{N^{k+1}} = \tilde m_k,\qquad \lim
_{N \to\infty} \frac{\mathbf E_N (p_k^2)}{N^{2(k+1)}} = \tilde m_k^2.
\]
\end{proposition}

The proof of this proposition is given in Section~\ref{proof43}.

\begin{lemma}
\label{geom}
There exists a constant $C_1>0$ such that $m_k < C_1^k$ for all $k \ge1$.
\end{lemma}

\begin{pf}
The general condition \eqref{cond-unit} implies that there exists a
constant $C>0$ such that
\[
t_k < C^k \qquad\mbox{for any $k \ge1$}.
\]
It is easy to see that if the coefficients of a formal power series are
majorated by a geometric progression, then the coefficients of the
inverse power series are also majorated by some geometric progression.
Therefore, the coefficients of the series
\[
v_0 (z) = \biggl( \frac{z}{1+z(1+z)(t_1 + t_2 z+ \cdots)} \biggr)^{(-1)}
\]
are majorated by a geometric progression. By definition, $m_k$ is the
coefficient of $z^{k+1}$ in $\log(1+v_0(z))$; this implies the
statement of the lemma.
\end{pf}

\begin{pf*}{Proof of Theorem~\ref{mainTh}}
Let $\lambda_i^{(N)}$, $i=1,2, \ldots, N$, be the coordinates of a
random signature distributed according to $M_N$. Note that Proposition~\ref{43} implies
%
\begin{eqnarray}
\label{conv-pk} \lim_{N \to\infty} \mathbf E_N \Biggl[
\frac{1}{N} \sum_{i=1}^N \bigl(
\lambda_i^{(N)} - i + 1/2 \bigr)^k \Biggr]&=&
m_k,
\\
\lim_{N \to
\infty} \mathbf E_N \Biggl[ \frac{1}{N}
\Biggl( \sum_{i=1}^N \bigl(
\lambda_i^{(N)} - i + 1/2 \bigr)^k \Biggr)
\Biggr]^2 &=& m_k^2.
\end{eqnarray}

It remains to apply Lemma~\ref{lemma-shod} (note that the existence of
the limit measure with moments $\{m_k\}$ follows from this lemma); the
conditions of the lemma hold due to~\eqref{conv-pk} and Lemma~\ref{geom}.
\end{pf*}

\subsection{Proof of Proposition \texorpdfstring{\protect\ref{41}}{4.1}}
\label{proof41}

Equations \eqref{difOp1} and \eqref{cond-unit} imply the following
formula for
$\mathbf E_N (p_k^{\#})$:
%
\begin{eqnarray}
\label{main-p} \mathbf E_N \bigl(p_k^{\#}\bigr)
&= &\sum_{1 \le i_1, i_2, \ldots, i_k \le N} \partial_{i_1 i_2}
\,\partial_{i_2 i_3} \cdots\partial_{i_k i_1}\nonumber\\
&&{}\times \exp \biggl( N \biggl(
\mathbf t_1(N) \operatorname{Tr}(X) + \frac{\mathbf t_2(N)}{2}
\operatorname{Tr}\bigl(X^2\bigr)
\\
&&\hspace*{36pt}\qquad{}+\cdots+ \frac{\mathbf t_r(N)}{r} \operatorname{Tr} \bigl(X^r\bigr) +
\cdots\biggr) \biggr)\bigg |_{X=0},
\nonumber
\end{eqnarray}
where the coefficients $\mathbf t_i (N)$ satisfy $\lim_{N \to\infty}
\mathbf
t_i (N) = t_i$ [the coefficients $t_i$ are given by \eqref{deriv}].

To deal with this formula, we need to introduce a bit of combinatorial
formalism. Below we use the term \emph{graph} to denote a finite connected
oriented graph, possibly with loops and multiple edges. A \textit{cycle}
in such
a graph is a closed oriented path without repeated edges. A cycle is
\textit{simple} if it does not contain repeated vertices. A cycle is said to
be \textit{Eulerian} if it contains all the edges of the graph. By an \textit{Eulerian
graph}, we mean a graph together with a distinguished enumeration of
the edges
such that it forms an Eulerian cycle (note that this slightly differs
from the
conventional terminology).

Let $\mathcal G_k$ denote the set of (equivalence classes of) Eulerian graphs
with $k$ edges. For $G\in\mathcal G_k$ we denote by $\bold e=(e_1,\ldots,e_k)$
the distinguished Eulerian cycle of $G$.

\begin{figure}

\includegraphics{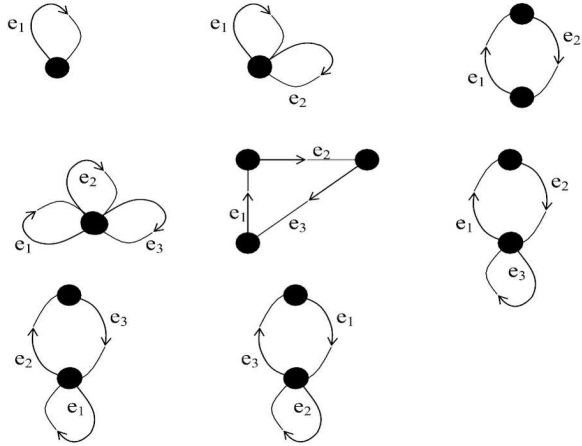}

\caption{All Eulerian graphs with $k=1,2,3$.}
\label{EulerGraphs}
\end{figure}

All Eulerian graphs with $k=1,2,3$ are shown in Figure~\ref{EulerGraphs}.

\begin{figure}

\includegraphics{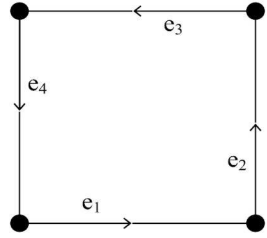}

\caption{A simple Eulerian graph with 4 vertices and 4 edges.}
\label{PolOrCycle}
\end{figure}

\begin{remark}\label{rem1}
There exists a one-to-one correspondence $G\leftrightarrow\pi$ between the
graphs $G\in\mathcal G_k$ and the set partitions of $[k]:=\{1,\ldots,k\}$.
Indeed, let us consider first the finest partition,
\[
\pi_0:=\{1\}\cup\{2\}\cup\cdots\cup\{k\}.
\]
By definition, the corresponding graph $G_0\leftrightarrow\pi_0$ is the
(unique) Eulerian graph with $k$ edges and $k$ vertices (see an example
in Figure~\ref{PolOrCycle}). Let
us enumerate the vertices of $G_0$ in such a way that
\[
e_1=(1\to2),\qquad\ldots,\qquad e_{k-1}=(k-1\to k),\qquad e_k=(k
\to1).
\]
Then, given an arbitrary set partition $\pi$ of $[k]$, we glue together the
vertices of $G_0$ corresponding to every block of $\pi$; the result is the
graph $G\leftrightarrow\pi$.

Equivalently, the vertices of $G$ are identified with the blocks of $\pi
$, and
the $i$th edge $e_i$ is directed from the block containing $i$ to that
containing $i+1$ (with the understanding that $k+1$ is identified with $1$).
\end{remark}

By $v(G)$, we will denote the number of vertices of $G$; this is the
same as the
number of blocks in the corresponding set partition $\pi$.

By a \textit{cycle structure} on $G$, we mean a partition $C=(C_1,\ldots
,C_p)$ of
the edge set $\{e_1,\ldots,e_k\}$ such that each block $C_j$ is a cycle;
we also
assume that the blocks are enumerated in the ascending order of their minimal
elements. Below we write the number of blocks by $p(C)$ and denote by $|C_j|$
the size of the $j$th block. The set of all cycle structures on $G$ is denoted
by $\mathcal C(G)$.

Note that cycle structures exist for every Eulerian graph $G$. For instance,
the Eulerian cycle $\bold e$ is itself a cycle structure with a single block.
Another example is obtained when one cuts $\bold e$ into simple cycles, which
is always possible, but sometimes can be made in different ways.

\begin{figure}[b]

\includegraphics{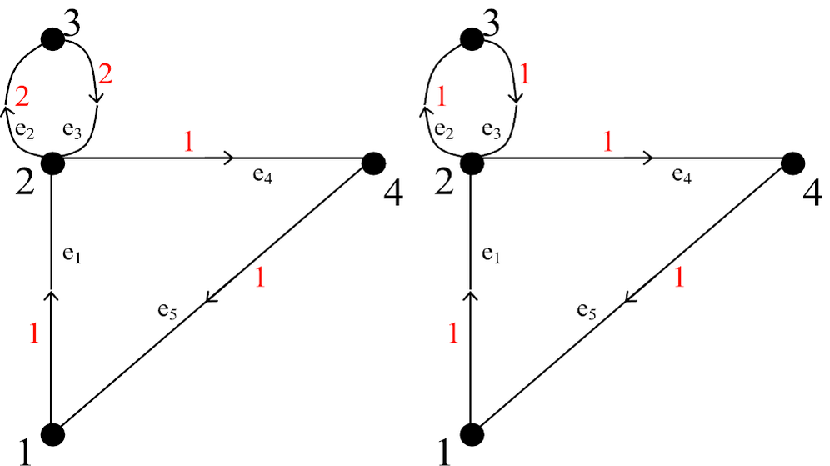}

\caption{All cycle structures on the Eulerian graph $1\to2\to3\to2\to
4\to1$.}
\label{OrCycle}
\end{figure}

Examples of cycle structures are shown in Figure~\ref{OrCycle}.

To shorten the notation, let us abbreviate
\[
\mathbf t_1:=\mathbf t_1(N),\qquad \mathbf t_2:=
\mathbf t_2(N), \ldots.
\]

\begin{lemma}\label{lemma1}
The right-hand side of \eqref{main-p} can be written in the form
%
\begin{equation}
\label{eq3} \sum_{G\in\mathcal G_k} N(N-1)\cdots\bigl(N-v(G)+1
\bigr)\sum_{C\in\mathcal
C(G)}\alpha(C)N^{p(C)}\prod
_{j=1}^{p(C)}\mathbf t_{|C_j|},
\end{equation}
$\alpha(C)$ is a coefficient depending on $C$ only. If all the blocks
of $C$
are simple cycles, then $\alpha(C)=1$.
\end{lemma}

\begin{pf}
\textit{Step} 1. Let us fix a sequence $\mathbf i = (i_1, i_2, \ldots,
i_k) \in[N]^k$.
The corresponding term in \eqref{main-p} can be written as
%
\begin{eqnarray}
\label{expr1} &&\partial_{i_1 i_2} \cdots\partial_{i_k i_1} \exp
\biggl( N \biggl( \mathbf t_1 \sum_{1 \le j_1 \le N}
x_{j_1 j_1} + \frac
{\mathbf
t_2}{2} \sum_{1 \le j_1,j_2 \le N}
x_{j_1 j_2} x_{j_2 j_1}
\nonumber
\\[-8pt]
\\[-8pt]
\nonumber
&&\hspace*{68pt}\qquad{}+ \cdots+ \frac{\mathbf t_r}{r} \sum_{1 \le j_1,j_2,
\ldots, j_r \le N}
x_{j_1 j_2} x_{j_2 j_3} \cdots x_{j_r j_1} + \cdots \biggr)
\biggr) \bigg|_{x=0} ,
\nonumber
\end{eqnarray}
where ``$x=0$'' means that finally all the $x$-variables are set to be
equal to
$0$.

The order of partial derivatives is not important; let us assume that one
applies $\partial_{i_1 i_2}$ first, then $\partial_{i_2 i_3}$, etc.
Since the
sum inside the exponential converges uniformly, we can differentiate this
expression term by term. Namely, each differentiation operator $\partial
$ can
be applied to one of the terms inside the exponential (as a result, a
pre-exponential polynomial appears), or it can be applied to a pre-exponential
factor which was brought down by previous differentiations. However,
due to the
final substitution $x=0$, a nonzero contribution can only come from
those terms
for which the pre-exponential factors do not contain the $x$-variables.

\textit{Step} 2. We will encode such terms by means of cycle structures.

First, we assign to $\mathbf i$ an Eulerian graph $G=G_{\mathbf i}$
with $k$ edges---the vertex set of $G$ is the subset of $[N]$ consisting of the numbers entering
the sequence $\mathbf i$, and the edges are
\begin{eqnarray*}
e_1&=&(i_1\to i_2),\qquad e_2=(i_2
\to i_3), \qquad\ldots,\qquad e_{k-1}=(i_{k-1}\to
i_k),\\
 e_k&=&(i_k\to i_1).
\end{eqnarray*}
In other words, we associate the edges with the $\partial$-operators in
\eqref{expr1}.

Next, given a term whose preexponential factor does not contain the
$x$-variables, we assign to it a partition $C=(C_1,\ldots,C_p)$ of the
edge set
$\{e_1,\ldots,e_k\}$ in the following way. The first block $C_1$ starts
with the
edge $e_1\leftrightarrow\partial_{i_1i_2}$, and the remaining edges correspond
to the $\partial$-operators killing the $x$-variables from the preexponential
factor that arises after application of $\partial_{i_1i_2}$ to the exponential.
The second block starts with the edge labeling the next $\partial$-operator
that is being applied to the exponential, etc.

We claim that $C$ is a cycle structure, that is, all blocks are cycles. Indeed,
a pre-exponential factor that may result from the application of
$\partial_{i_1i_2}$ to the exponential always has the form
%
\begin{equation}
\label{eq1} \frac{N}r \mathbf t_r x_{j_1j_2}
\cdots\widehat{x_{j_mj_{m+1}}}\cdots x_{j_rj_1}, \qquad j_m=i_1,
j_{m+1}=i_2
\end{equation}
(with the understanding that $m+1=1$ if $m=r$; the hat over $x_{i_mi_{m+1}}$
means that this variable has to be omitted). Our assumption is that $r=|C_1|$
and the $\partial$-operators corresponding to the edges from $C_1$ different
from $e_1$ kill all the $x$-variables from the above monomial. But this just
means that the edges of $C_1$ form a cycle.

For the blocks $C_2$, $C_3$, etc. the argument is the same.

\textit{Step} 3. The reasoning of step 2 shows that the quantity \eqref
{expr1} can be
represented as the sum of contributions coming from various cycle structures
$C\in\mathcal C(G)$. Let us fix $C=(C_1,\ldots,C_p)$ and analyze its
contribution in more detail. Assume first that all the cycles are
simple. Let
us focus on the first cycle and keep the notation of step~2. The fact that
$C_1$ is simple just means that the indices $j_1,\ldots,j_r$ must be pairwise
distinct. Therefore, given a simple cycle $C_1$, there are exactly $r=|C_1|$
eligible $r$-tuples $(j_1,\ldots,j_r)$ that correspond to values
$m=1,\ldots,r$.
Then the summation over these $r$ variants results in the cancellation
of the
factor $r$ in the denominator of \eqref{eq1}. The same argument
applies to all the cycles, and we finally obtain that the whole
contribution of
$C$ is equal to
\[
N^p\mathbf t_{|C_1|}\cdots\mathbf t_{|C_p|},
\]
as desired.

In the general case, when the cycles are not necessarily simple, we
argue as
above, and the only difference is that the contribution of $C$ may
involve a
constant numeric factor $\alpha(C)$. For instance, if the graph $G$ has
a single
vertex and $k$ loops, then there is a single one-component cycle structure
whose contribution equals $(k-1)!N\mathbf t_k$, so that in this case
$\alpha(C)=(k-1)!$.

\textit{Step} 4. We have explained the origin of the interior sum in
\eqref{eq3}. It
remains to explain the exterior sum, and this is easy. Namely, we
observe that
the whole contribution of a given $k$-tuple $\mathbf i\in[N]^k$ depends solely
on the equivalence class of the corresponding Eulerian graph
$G_{\mathbf i}$.
Indeed, two $k$-tuples producing equivalent graphs can be transformed
to each
over by a permutation of $[N]$, which does not affect the quantity
\eqref{expr1}. Finally, given $G\in\mathcal G_k$, the number of $k$-tuples
$\mathbf i\in[N]^k$ such that $G_{\mathbf i}$ is equivalent to $G$ is
equal to
\[
N(N-1)\cdots\bigl(N-v(G)+1\bigr)
\]
(to see this one may use Remark~\ref{rem1}). This completes the proof.
\end{pf}

Let us rewrite \eqref{eq3} as
\[
\sum_{(G,C)} \alpha(C) N(N-1)\cdots\bigl(N-v(G)+1
\bigr)N^{p(C)}\prod_{j=1}^{p(C)}
\mathbf t_{|C_j|},
\]
where the summation is taken over all pairs $(G,C)$ such that $G\in
\mathcal
G_k$ and $C\in\mathcal C(G)$. For $N$ large, the contribution of a
fixed pair
$(G,C)$ grows as $N^{v(G)+p(C)}$, so that the leading part in the asymptotics
comes from the pairs with the maximal possible value of the quantity
$v(G)+p(C)$. Our goal now is to describe such pairs.

Assume $A$ is an ordered set and $A_1\subset A$ and $A_2\subset A$ are two
nonempty disjoint subsets. Then $A_1$ and $A_2$ are said to be \textit{crossing}
if there exists a quadruple $a<b<c<d$ of elements such that $a$ and $c$
are in
one of these subsets while $b$ and $d$ are in another subset; otherwise $A_1$
and $A_2$ are said to be \textit{noncrossing}. Next, a \textit{noncrossing
partition} of $A$ is a set partition of $A$ whose blocks are pairwise
noncrossing.

By the very definition, every cycle structure $C=(C_1,\ldots,C_p)$ on a graph
$G\in\mathcal G_k$ is a partition of the set $\{e_1,\ldots,e_k\}$. We introduce
the natural order $e_1<\cdots<e_k$ on the Eulerian cycle, so that
$\{e_1,\ldots,e_k\}$ becomes an ordered set isomorphic to~$[k]$.

\begin{lemma}\label{lemma2}
Let us fix $k=1,2,\ldots$ and let $(G,C)$ range over the set of pairs
such that
$G\in\mathcal G_k$ and $C\in\mathcal C(G)$.

Then the maximal possible value of the quantity $v(G)+p(C)$ is equal to $k+1$.
It is attained exactly for those pairs $(G,C)$ for which all the cycles
of $C$
are simple and the set partition $\sigma(C)$ is noncrossing.

Moreover, under the identification $\{e_1,\ldots,e_k\}\leftrightarrow
[k]$ of
ordered sets, for every noncrossing partition $\sigma$ of the set
$[k]$, there
exists exactly one pair $(G,C)$ such that $v(G)+p(C)=k+1$ and
$C\leftrightarrow\sigma$.
\end{lemma}

\begin{pf}
\textit{Step} 1. Let us fix a pair $(G,C)$ with $C=(C_1,\ldots,C_p)$,
and estimate
$v(G)+p(C)$.

Let us observe that $C_2$ always has a common vertex with $C_1$, $C_3$
has a
common vertex with $C_1\cup C_2$, and so on. Indeed, this follows from
the very
definition of a cycle structure (in particular, we use the fact the
cycles in
$C$ are enumerated in the ascending order of their minimal elements).

Let $v( \cdot )$ stand for the number of vertices in a given cycle or
a union
of cycles. We have
\[
v(C_1)\le|C_1|,\qquad \ldots,\qquad v(C_p)
\le|C_p|
\]
and, by virtue of the above observation,
\begin{eqnarray*}
v(C_1\cup\cdots\cup C_m)&\le& v(C_1\cup
\cdots\cup C_{m-1})+v(C_m)-1\\
&\le& v(C_1\cup
\cdots\cup C_{m-1})+|C_m|-1
\end{eqnarray*}
for $m=2,\ldots,p$. Since $|C_1|+\cdots+|C_p|=k$, it follows that
\[
v(G)=v(C_1\cup\cdots\cup C_p)\le k-(p-1),
\]
so that $v(G)+p\le k+1$.

Moreover, the equality $v(G)+p=k+1$ is attained if and only if the following
two conditions are satisfied:
\begin{longlist}[(1)]
\item[(1)] $v(C_m)=|C_m|$ for every $m=1,\ldots,p$, which is equivalent to
saying that all cycles are simple.

\item[(2)] For every $m=2,\ldots,p$, the cycle $C_m$ has a single
common vertex
with the union $C_1\cup\cdots\cup C_{m-1}$.
\end{longlist}

\textit{Step} 2. Let us assume that $(G,C)$ is such that $C$ satisfies
condition (1)
above; we are going to show that $C$ satisfies condition (2) if and
only if $C$
is noncrossing.

The key observation is that if $(G,C)$ is such that $C$ satisfies both
(1) and
(2), then removing the last cycle $C_p$ we still get a pair $(G',C')$
with the
same properties. Likewise, if $\sigma$ is a noncrossing set partition, then
removing its last block we still get a noncrossing partition $\sigma'$ (we
always assume that the blocks are ordered according to the order of their
minimal elements).

This suggests the idea to prove the desired claim by induction on $p$, the
number of blocks. The base of induction is obvious: if $p=1$, then
there is
nothing to prove. To justify the induction step, we observe that the possible
transitions $(G',C')\to(G,C)$ preserving property (2) are directed by exactly
the same mechanism as the possible transitions $\sigma'\to\sigma$ preserving
the noncrossing property.

Indeed, in the first case, we may insert a simple cycle of length
$|C_p|$ at
any place of the Eulerian cycle of $G'$ which is after the minimal edge of
$C_{p-1}$ (which is the last cycle of $C'$). Likewise, in the second
case, we
may insert a block of the same size after the minimal element of the
last block
of $\sigma'$. (Let us emphasize that in both cases, we have to insert a new
cycle/block as a whole.)

\textit{Step} 3. The argument of step 3 shows that both the pairs
$(G,C)$ satisfying
conditions (1) and (2), and the noncrossing set partitions $\sigma$ can be
obtained by one and the same recursive procedure. This completes the proof.
\end{pf}

\begin{figure}

\includegraphics{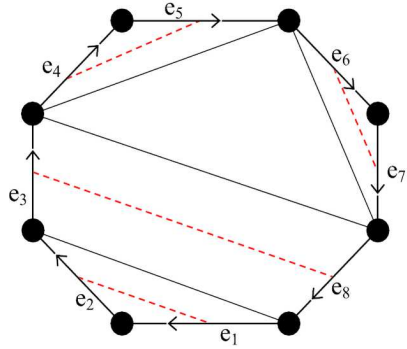}

\caption{A noncrossing partition of edges (dashed lines) gives rise to
a noncrossing partition of vertices (solid lines).}
\label{Krew}
\end{figure}

\begin{remark}
The recursive procedure described above assigns a pair $(G,C)$ to every
noncrossing partition $\sigma$ of the set $[k]$. On the other hand, according
to Remark~\ref{rem1}, the graph $G$ is completely determined by a set partition
$\pi$ of $[k]$. One can show that the correspondence $\sigma\mapsto\pi$ that
arises in this way is just the \textit{complementation} operation first
discovered by Kreweras \cite{Kreweras}: it is a nontrivial involution
on the
set of noncrossing partitions of $[k]$ (see an example in Figure~\ref{Krew}).
\end{remark}

Denote by $\operatorname{NC}_k$ the set of noncrossing partitions of $[k]$.
Define the \textit{weight} of a partition
$\sigma=(\sigma_1\cup\cdots\cup\sigma_p)\in\operatorname{NC}_k$ as the monomial
\[
\operatorname{wt}(\sigma):=t_{|\sigma_1|}\cdots t_{|\sigma_p|}.
\]
Lemmas \ref{lemma1} and \ref{lemma2} show that the leading term of the
large-$N$ asymptotics can be written as
%
\begin{equation}
N^{k+1}\sum_{\sigma\in\operatorname{NC}_k}\operatorname{wt}(
\sigma).
\end{equation}

\begin{lemma}
For any $k \ge1$, we have
\[
\sum_{\sigma\in\operatorname{NC}_k}\operatorname{wt}(\sigma)=
\frac{1}{k+1} \bigl[u^k\bigr] \bigl\{ \bigl(1 + t_1
u + t_2 u^2 + \cdots\bigr)^{k+1} \bigr\}.
\]
\end{lemma}

\begin{pf}
Given $\sigma\in\operatorname{NC}_k$, let $(1^{s_1}2^{s_2}\cdots)$
denote the
corresponding ordinary partition of the \textit{number} $k$, written in the
multiplicative notation; this means that $\sigma$ has exactly $s_i$
blocks of
size $i$, where $i=1,2,\ldots.$ We say that $(1^{s_1}2^{s_2}\cdots
k^{s_k})$ is
the \textit{type} of $\sigma$. Obviously,
\[
\operatorname{wt}(\sigma)=t_1^{s_1} t_2^{s_2}
\cdots t_k^{s_k}.
\]
Therefore, we have to prove that
\[
\sum_{\sigma\in\operatorname{NC}_k} t_1^{s_1}
t_2^{s_2}\cdots t_k^{s_k}=
\frac{1}{k+1} \bigl[u^k\bigr] \bigl\{ \bigl(1 + t_1
u + t_2 u^2 + \cdots\bigr)^{k+1} \bigr\}.
\]

Now we apply Exercise 5.35a in Stanley \cite{Sta}, which says that the number
of partitions $\sigma\in\operatorname{NC}_k$ of a given type
$(1^{s_1}2^{s_2}\cdots k^{s_k})$ is equal to
\begin{eqnarray}
\frac{k(k-1)\cdots(k-\ell+2)}{s_1!s_2!\cdots s_k!} =
\frac{1}{k+1}\frac{(k+1)k(k-1)\cdots(k-\ell+2)}{s_1!s_2!\cdots}, \nonumber\\
\eqntext{\ell:=s_1+s_2+
\cdots+s_k.}
\end{eqnarray}
This is equivalent to the desired formula.
\end{pf}

\subsection{Proof of Proposition \texorpdfstring{\protect\ref{42}}{4.2}}
\label{proof42}
Let us abbreviate $l:=l(\rho)$. Relations \eqref{difOp1} and~\eqref{cond-unit}
imply
%
\begin{eqnarray}
\label{main-p-rho}\lim_{N \to\infty} \frac{\mathbf E_N
(p_{\rho}^{\#})}{N^{\operatorname{wt}(\rho)}}& =& \lim
_{N \to\infty} \frac
{1}{N^{k_1+k_2 +
\cdots+ k_{l} + l}} \nonumber\\
&&\times{}\sum_{\mathbf i \in[N]^k}
\partial_{i_1 i_2}\, \partial_{i_2 i_3} \cdots\partial_{i_{k_1} i_1}
\,\partial_{i_{k_1+1} i_{k_1+2}} \cdots
\nonumber
\\[-8pt]
\\[-8pt]
\nonumber
&&{}\times\exp\biggl( \mathbf t_1(N) \biggl( \sum
_{j} x_j \biggr) + \frac{\mathbf t_2(N)}{2} \biggl(
\sum_{j_1, j_2} x_{j_1 j_2} x_{j_2
j_1}
\biggr) + \cdots
\\
&&\hspace*{46pt}{}+ \frac{\mathbf t_r(N)}{r} \biggl( \sum_{j_1, j_2, \ldots, j_r}
x_{j_1 j_2} x_{j_2 j_3} \cdots x_{j_r j_1} \biggr) + \cdots
\biggr) \bigg|_{x_{ij} \equiv0},
\nonumber
\end{eqnarray}
where the coefficients $\mathbf t_i (N)$ satisfy $\lim_{N \to\infty}
\mathbf t_i (N) = t_i$ [the numbers $t_i$ were defined in \eqref{deriv}].

We shall deal with this formula in the same way as in Section~\ref{proof41}. To
every sequence $\mathbf i=(i_1,\ldots,i_k)\in[N]^k$ we assign an
oriented graph
$G_{\mathbf i}$ whose edges correspond to the $\partial$-operators from
\eqref{main-p-rho}. This graph is composed from $l$ Eulerian graphs,
which may
be glued together or disjoint, depending on whether the subsequences
%
\begin{equation}\qquad
\label{eq2} (i_1,\ldots,i_{k_1}),\hspace*{-1pt}\qquad (i_{k_1+1},\hspace*{-1pt}
\ldots,i_{k_1+k_2}), \hspace*{-1pt}\qquad\ldots,\qquad (i_{k_1+\cdots+
k_{l-1}+1}, \ldots,i_k)\hspace*{-10pt}
\end{equation}
have common indices or not.

First, let us consider the case when there are no common indices, so
that the
corresponding Eulerian graphs are pairwise disjoint. Then the differential
operators from different graphs are applied to nonintersecting sets of
$x$-variables, and the arguments of Section~\ref{proof41} show that the total
contribution from such $\mathbf i$'s equals
\[
c_{k_1} c_{k_2} \cdots c_{k_{l}} N^{(k_1+1) + (k_2+1) + \cdots+(k_l+1)} +
O\bigl(N^{k+l-1}\bigr).
\]

It remains to show that the contribution from the remaining sequences
$\mathbf
i$ [those for which the subsequences in \eqref{eq2} have common
indices] has
lower degree in~$N$.

To simplify the argument, let us assume that $l=2$, so that $k=k_1+k_2=r+(k-r)$.
Thus, there are two subsequences in \eqref{eq2}, which we denote as
\[
(i_1,\ldots,i_r),\qquad (i_{r+1},\ldots,
i_k),
\]
and these two subsequences share a common index, say $i_a=i_{r+b}$ for some
$a\in\{1,\ldots,r\}$ and $b\in\{1,\ldots,k-r\}$.

Then it is readily seen that the term corresponding to the differential
operator
\[
\partial_{i_1 i_2} \,\partial_{i_2 i_3} \cdots\partial_{i_r i_1}
\,\partial_{i_{r+1} i_{r+2}} \,\partial_{i_{r+2} i_{r+3}} \cdots\partial _{i_k i_{r+1}}
\]
is equal to the contribution of a \emph{single} Eulerian graph with $k$ edges,
corresponding to the sequence
\[
i_1,\ldots,i_a, i_{r+b+1}, i_{r+b+2},
\ldots, i_k, i_{r+1}, i_{r+2}, \ldots,
i_{r+b}, i_{a+1}, i_{a+2}, \ldots,
i_r, i_1.
\]
Therefore, this contribution has order at most $N^{k+1}$, which is less than
$N^{k+l}=N^{k+2}$.

The same argument holds when $l>2$ as well.

\subsection{Proof of Proposition \texorpdfstring{\protect\ref{43}}{4.3}}
\label{proof43}

Recall that we consider the functions $p_k$ as random variables on the
probability space $(\mathbb{GT}_N, M_N)$. First, let us prove that
after scaling the functions $p_k$ converge to constants in $L^2$.

\begin{lemma}
There exist constants $\bar m_k$, $k =1,2, \ldots,$ such that for any $k
\ge1$
\[
\lim_{N \to\infty} \frac{\mathbf E_N (p_k)}{N^{k+1}} = \bar m_k,\qquad \lim
_{N \to\infty} \frac{\mathbf E_N (p_k^2)}{N^{2(k+1)}} = \bar m_k^2.
\]
\end{lemma}

\begin{pf}
Let $f \in\operatorname{Sym}^*$ be arbitrary. Since $f$ is a linear
combination of $p_{\rho}^{\#}$'s with $\operatorname{wt}(\rho) \le \operatorname{wt}(f)$, and there exist
limits of $\mathbf E_N (p_{\rho}^{\#})/N^{k+1}$, we obtain
\[
\lim_{N \to\infty} \frac{\mathbf E_N (f)}{N^{k+1}} = a_f,
\]
for some constants $a_f$.

It is known (see \cite{IvaOls}) that
\[
p_{\rho_1}^{\#} p_{\rho_2}^{\#} =
p_{\rho_1 \cup\rho_2}^{\#} + \mbox {lower weight terms},
\]
where $\rho_1 \cup\rho_2$ stands for the union of the partitions $\rho
_1$ and $\rho_2$, and ``lower weight terms'' denotes terms with weight
$\le \operatorname{wt}(p_{\rho_1}^{\#}) + \operatorname{wt}(p_{\rho_2}^{\#}) -1$.
Hence,
\[
\lim_{N \to\infty} \frac{\mathbf E_N  ( p_{\rho_1}^{\#} p_{\rho
_2}^{\#}  )}{N^{\operatorname{wt}(p_{\rho_1}^{\#}) + \operatorname{wt}(p_{\rho_2}^{\#})}} = \lim_{N \to\infty}
\frac{\mathbf E_N  ( p_{\rho_1 \cup\rho_2}^{\#}
 )}{N^{\operatorname{wt}(p_{\rho_1}^{\#}) + \operatorname{wt}(p_{\rho_2}^{\#})}}.
\]

This equality and Proposition~\ref{42} imply that
\[
\lim_{N \to\infty} \frac{\mathbf E_N (f^2)}{N^{2 \operatorname{wt}(f)}} = a_f^2.
\]
Therefore, the functions $f$ converge to $a_f$ in $L^2$.

Choosing the function $p_k$ as $f$ we obtain the statement of the
lemma. By $\bar m_k$ we denote the limit constant.
\end{pf}

It remains to prove that $\bar m_k = \tilde m_k$ for all $k= 1,2, \ldots
$ (recall that the constants $\tilde m_k$ were defined in Section~\ref{statement}).

Consider formal power series of the form
\[
a(z) = a_1 z + a_2 z^2 + \cdots,\qquad
a_i \in\mathbb R.
\]

Recall that a series of this form is invertible if and only if $a_1 \ne
0$. Let $a^{(-1)}(z)$ denote the inverse of $a(z)$, that is, $a^{(-1)}
(a(z)) = z$.
Set
\[
\bar A (z) = \bar m_1 z^2 + \bar m_2
z^3 + \bar m_3 z^4 + \cdots
\]
and
\[
C(z) = 1 + c_1 z^2 + c_2 z^3 +
c_3 z^4 + \cdots.
\]

\begin{lemma}
The formal power series $z \exp(\bar A(z))$ and $z/C(z)$ are inverse
to each other.
\end{lemma}

\begin{pf}
Recall that [see \eqref{change}]
\[
\mathbf p_k = \frac{1}{k+1} \bigl[ u^{k+1}\bigr] \bigl
\{ \bigl(1 + \mathbf p_1^{\#} u^2 + \mathbf
p_2^{\#} u^3 + \cdots\bigr)^{k+1}
\bigr\} + \mbox{lower weight terms},
\]
where ``lower weight terms'' denotes terms with weight $\le k$. Since
\[
\mathbf E_N (f) = O\bigl(N^{\operatorname{wt}(f)}\bigr),\qquad f \in\mathbb A(N),
\]
the ``lower weight terms'' do not affect the asymptotics of $\mathbf
E_N (p_k)$, and we have
\[
\bar m_k = \frac{1}{k+1} \bigl[ u^{k+1}\bigr] \bigl\{
\bigl(1 + c_1 u^2 + c_2 u^3 +
\cdots \bigr)^{k+1} \bigr\}.
\]
The lemma follows from this formula; cf. \cite{IvaOls}, Propositions 3.6, 3.7.
\end{pf}

Let us find an expression for $C(z)$ using the formula for $c_k$'s
given by Proposition~\ref{41}.

\begin{lemma}
We have
\[
C(z) = 1 - z + \biggl( \frac{z}{1+t_1 z+ t_2 z^2 + \cdots} \biggr)^{(-1)}.
\]
\end{lemma}

\begin{pf}
This is a direct consequence of the Lagrange inversion formula (see,
e.g., \cite{Sta}, Theorem~5.4.2).
\end{pf}

Two previous lemmas imply that
\[
\bar m_1 z^2 + \bar m_2 z^3 +
\cdots= \log \biggl( \frac{1}{z} \biggl( \frac{z}{1 - z +  ( {z}/{(1+t_1 z+ t_2 z^2 + \cdots)}
)^{(-1)} }
\biggr)^{(-1)} \biggr).
\]

In order to show that $\bar m_k=\tilde m_k$, $k=1,2,\ldots,$ we prove
the equality of their generating functions. Formulas \eqref{Sz} and
\eqref{Stilt} imply
\[
\tilde m_1 z^2 + \tilde m_2 z^3
+ \cdots= \log \biggl( \frac{z}{1 +
z(1+z)(t_1 + t_2 z+ \cdots)} \biggr)^{(-1)} - \log(1+z).
\]

Therefore, the following lemma completes the proof.

\begin{lemma}
We have
%
\begin{eqnarray}
\label{eq-series} &&\biggl( \frac{z}{1 + z(1+z)(t_1 + t_2 z+ \cdots)}
\biggr)^{(-1)} +1
\nonumber
\\[-8pt]
\\[-8pt]
\nonumber
&&\qquad =
\frac{z+1}{z} \biggl( \biggl( \frac{z}{1 - z +  ( {z}/{(1+t_1 z+
t_2 z^2 + \cdots)}  )^{(-1)} } \biggr)^{(-1)}
\biggr).
\end{eqnarray}
\end{lemma}

\begin{pf}
It is easy to see that both series have the form
\[
a_0 + a_1 z + a_2 z^2 +
\cdots,
\]
and the coefficients $a_0$ and $a_1$ of both of the series are equal to
1. Let us prove the equality of the coefficients of $z^n$, $n \ge2$.

Recall that the Lagrange inversion formula has the following form (see,
e.g.,~\cite{Sta}, Theorem~5.4.2)
\[
n \bigl[z^n\bigr] F^{(-1)} (z)^k = k
\bigl[z^{n-k}\bigr] \biggl( \frac{z}{F(z)} \biggr)^n.
\]
Let $s(z)$ denote the series $1+ t_1 z + t_2 z^2 + \cdots.$ We have
\[
1 + z(1+z) (t_1 + t_2 z + \cdots) = s(z) + z s(z) - z.
\]
Hence, the coefficient of $z^n$ in the left-hand side of \eqref
{eq-series} can be written in the following form:
%
\begin{eqnarray}
\label{l-side} &&\frac{1}{n} \bigl[z^{n-1}\bigr] \bigl(s(z)+z s(z)
- z\bigr)^n \nonumber\\
&&\qquad= \frac{1}{n} \sum_{i_1+i_2+i_3 = n; i_1 \ge1}
\frac{n!}{i_1 ! i_2 ! i_3 !} \bigl[z^{n-1}\bigr] s(z)^{i_1}
z^{i_2} s(z)^{i_2} (-1)^{i_3} z^{i_3}
\\
&&\qquad= \sum_{i_1+i_2+i_3=n; i_1 \ge1} \frac{(n-1)!}{ i_1 ! i_2 ! i_3 !}
\bigl[z^{i_1-1}\bigr] s(z)^{i_1 + i_2} (-1)^{i_3}.
\nonumber
\end{eqnarray}

Now let us consider the expression in the right-hand side of $\eqref
{eq-series}$. The Lagrange inversion formula implies that the
coefficient of $z^n$ in the right-hand side can be written as
\[
\frac{1}{n} \bigl[z^{n-1}\bigr] C(z)^n +
\frac{1}{n+1} \bigl[z^n\bigr] C(z)^{n+1}.
\]
Since
\[
C(z) = 1 - z + \biggl( \frac{z}{s(z)} \biggr)^{(-1)},
\]
the coefficient of $z^l$ (for any $l \ge1$) in an arbitrary power of
$C(z)$ can be found by the Lagrange inversion formula again. We obtain
\[
\bigl[z^l\bigr] \biggl( \biggl( \frac{z}{s(z)}
\biggr)^{(-1)} (z) \biggr)^k = \frac
{k}{l}
\bigl[z^{k-l}\bigr] s(z)^l.
\]
Therefore,
\begin{eqnarray*}
&&\frac{1}{n} \bigl[z^{n-1}\bigr] C(z)^n\\
&&\qquad =
\frac{1}{n} \biggl( \sum_{i_1+i_2+i_3 =
n; i_1 \ge1, i_3 \ge1}
(-1)^{i_2} \frac{i_3}{n-1-i_2} \bigl[z^{i_1-1}\bigr]
s(z)^{n-1-i_2} \frac{n!}{i_1 ! i_2 ! i_3 !}
\\
&&\hspace*{265pt}{}+ (-1)^{n-1} n \biggr).
\end{eqnarray*}
This formula with $n$ replaced by $n+1$ reads
\begin{eqnarray*}
&&\frac{1}{n+1} \bigl[z^{n}\bigr] C(z)^{n+1}
\\
&&\qquad=
\frac{1}{n} \biggl( \sum_{j_1+j_2+j_3
= n+1; j_1 \ge1, j_3 \ge1}
(-1)^{j_2} \frac{j_3}{n-j_2} \bigl[z^{j_1-1}\bigr]
s(z)^{n-j_2} \frac{(n+1)!}{j_1 ! j_2 ! j_3 !}
\\
&&\hspace*{222pt}\qquad{}+ (-1)^{n} (n+1) \biggr).
\end{eqnarray*}

Now add the two equalities above and combine the coefficients of $[z^a]
s(z)^b$ for all $a$ and $b$. We obtain the sum
\[
\sum_{j_1 + j_2 + j_3 =n; j_1 \ge0} (-1)^{j_2} \frac{(n-1)!}{j_1 ! j_2
! j_3 !}
\bigl[z^{j_1-1}\bigr] s(z)^{j_1 + j_3}.
\]
Changing the notation of indices in the summation, we see that this
expression coincides with \eqref{l-side}. This completes the proof.
\end{pf}

\begin{appendix}
\section*{Appendix: Examples of limit shapes}\label{app}
In this section, we consider several examples of sequences $\omega
=\omega(N)$ which satisfy the main condition \eqref{cond-unit}.

(a) \textit{One-sided Plancherel character}.
Let $\gamma^+ = \gamma N$, where $\gamma$ is a fixed constant, and all
other Voiculescu's parameters are equal to 0. In this case, the main
condition holds with $t_1 = \gamma$ and $t_k =0$, for $k \ge2$. Then
we have
\[
Q(z) = 1+ \gamma z(1+z).
\]

It follows that
\[
v_0 (z) = \frac{1 - \gamma z - \sqrt{y^2 (z^2 - 4 \gamma) - 2 \gamma z
+1}}{2 \gamma z}.
\]

Using \eqref{Sz}, one can derive the expression for the Stieltjes transform:
\[
S \biggl( \frac{1}{z} \biggr) = \mathrm{Stil}^{\mathrm{Planch}} (z) = \log
\frac{z+\gamma
- \sqrt{(z-\gamma)^2 - 4 \gamma}}{ 2 \gamma}.
\]

Given the Stieltjes transform of a measure, there is a standard way to
compute the density of this measure; see, for example, \cite{AGZ}, Section~2.4, and the end of Section~\ref{DetPr} above. After
computations, we obtain that for $\gamma>1$ we have
\[
d^{\mathrm{Planch}} (x) = \frac{1}{\pi} \arccos\frac{x+\gamma}{2 \sqrt{\gamma(x+1)}}\qquad \mbox{for
$x \in[\gamma- 2 \sqrt{\gamma}; \gamma+ 2\sqrt{\gamma} ]$ },
\]
and for $\gamma< 1$ we have
\[
d^{\mathrm{Planch}} (x) = %
\cases{\displaystyle\frac{1}{\pi} \arccos
\frac{x+\gamma}{2 \sqrt{\gamma(x+1)}}, & \quad$\mbox{for $x \in[\gamma- 2 \sqrt{\gamma}; \gamma+ 2\sqrt{
\gamma} ]$,} $\vspace*{2pt}
\cr
1, &\quad $\mbox{for $x \in[-1; \gamma- 2 \sqrt{
\gamma}]$}$.} %
\]

\begin{figure}

\includegraphics{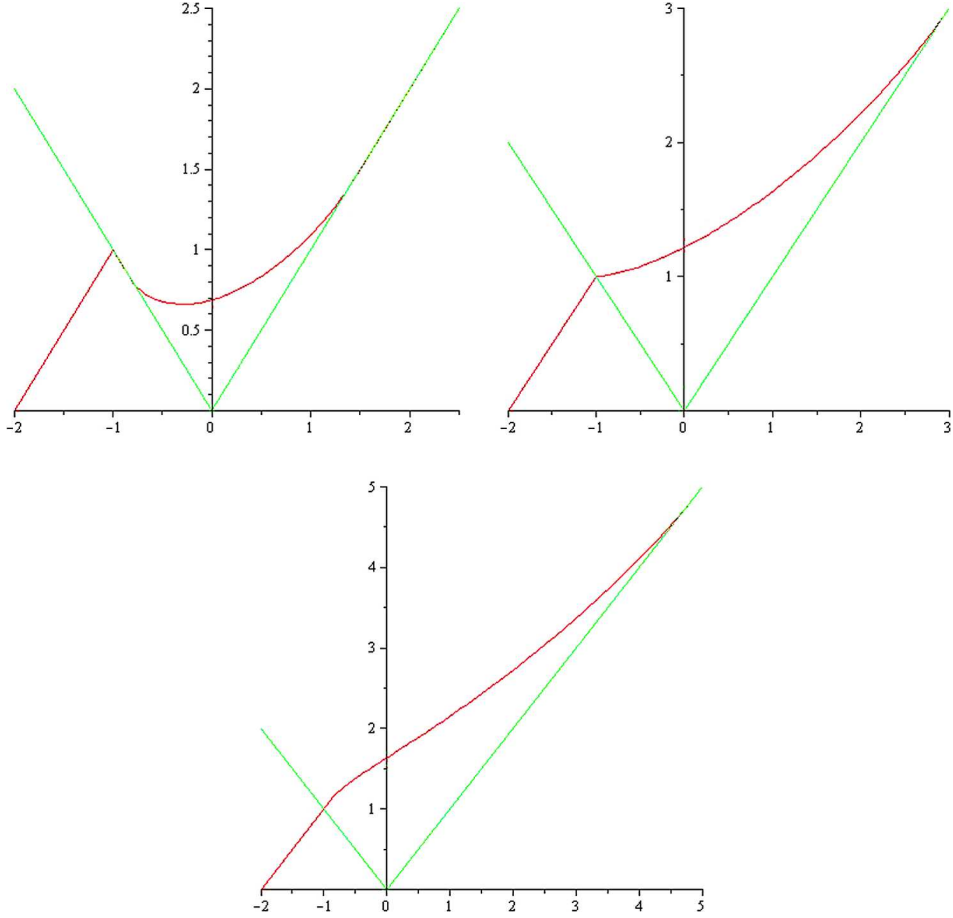}

\caption{Limit shapes for the one-sided Plancherel character with
$\gamma=0.3, 1,2$, respectively.}\label{limitshapes1P}
\end{figure}

Examples of these limit shapes are shown in Figure~\ref{limitshapes1P}.

After rescaling, these limit shapes coincide with Biane's limit shapes
(see \cite{Bia01}).

(b) \textit{One multiple $\alpha^+$-parameter}.
Assume that $\alpha_1^+ = \alpha_2^+ = \cdots= \alpha_{[a N]}^+ = \alpha
$, and all other Voiculescu's parameters are equal to 0. Note that we
fix two different real numbers $a$ and $\alpha$. Then $t_1 = a \alpha$,
$t_2 = a \alpha^2, \ldots,t_k = a \alpha^k, \ldots.$

In this case, we have
\[
Q(z) = 1 + \frac{z(1+z) a \alpha}{1 - \alpha z}.
\]

After computations, we obtain
\[
\mathrm{Stil}^{\mathrm{multi}\mbox{-}\alpha} (z) = \log\frac{\alpha(a+1) + (2 \alpha+1)z -
\sqrt{(z-\alpha(a+1))^2 - 4 a \alpha(\alpha+1)}}{2 \alpha(a+z)}.
\]

The limiting density is given by the following formulas.

For $a \ge(\alpha+1)/\alpha$, we have
\begin{eqnarray}
d^{\mathrm{multi}\mbox{-}\alpha} (x) = \frac{1}{\pi} \arccos\frac{\alpha(a+1) + (2
\alpha+1) x}{2 \sqrt{\alpha(\alpha+1) (x+1) (x+a)}},
\nonumber\\
\eqntext{x \in\bigl[\alpha(a+1) - 2 \sqrt{a \alpha(\alpha+1)}; \alpha(a+1) + 2
\sqrt{a \alpha(\alpha+1)}\bigr]. }
\end{eqnarray}

For $\alpha/ (\alpha+1) \le a \le(\alpha+1)/\alpha$, we have
\begin{eqnarray*}
&&d^{\mathrm{multi}\mbox{-}\alpha} (x)
\nonumber
\\[-8pt]
\\[-8pt]
\nonumber
&&\qquad= %
\cases{ 1, \qquad\mbox{$x \in\bigl[-1; \alpha(a+1) -
2 \sqrt{a \alpha(\alpha +1)}\bigr]$}, \vspace*{2pt}
\cr
\displaystyle\frac{1}{\pi} \arccos
\frac{\alpha(a+1) + (2 \alpha+1) x}{2 \sqrt
{\alpha(\alpha+1) (x+1) (x+a)}}, \vspace*{2pt}
\cr
\qquad\quad x \in\bigl[\alpha(a+1) - 2 \sqrt{a \alpha(
\alpha+1)}; \alpha(a+1) + 2 \sqrt{a \alpha(\alpha +1)}\bigr]. } %
\end{eqnarray*}
Finally, for $a \le\alpha/ (\alpha+1)$ we have
\begin{eqnarray*}
&&d^{\mathrm{multi}\mbox{-}\alpha} (x) \\
&&\qquad= %
\cases{ 1, \qquad\mbox{$x \in[-1; -a]$},
\vspace*{2pt}
\cr
\displaystyle\frac{1}{\pi} \arccos\frac{\alpha(a+1) + (2 \alpha+1) x}{2 \sqrt
{\alpha(\alpha+1) (x+1) (x+a)}}, \vspace*{2pt}
\cr
\qquad\quad\mbox{$x \in\bigl[\alpha(a+1) - 2 \sqrt{a \alpha(\alpha+1)}; \alpha(a+1) + 2
\sqrt{a \alpha(\alpha +1)}\bigr]$. }} %
\end{eqnarray*}

\begin{figure}

\includegraphics{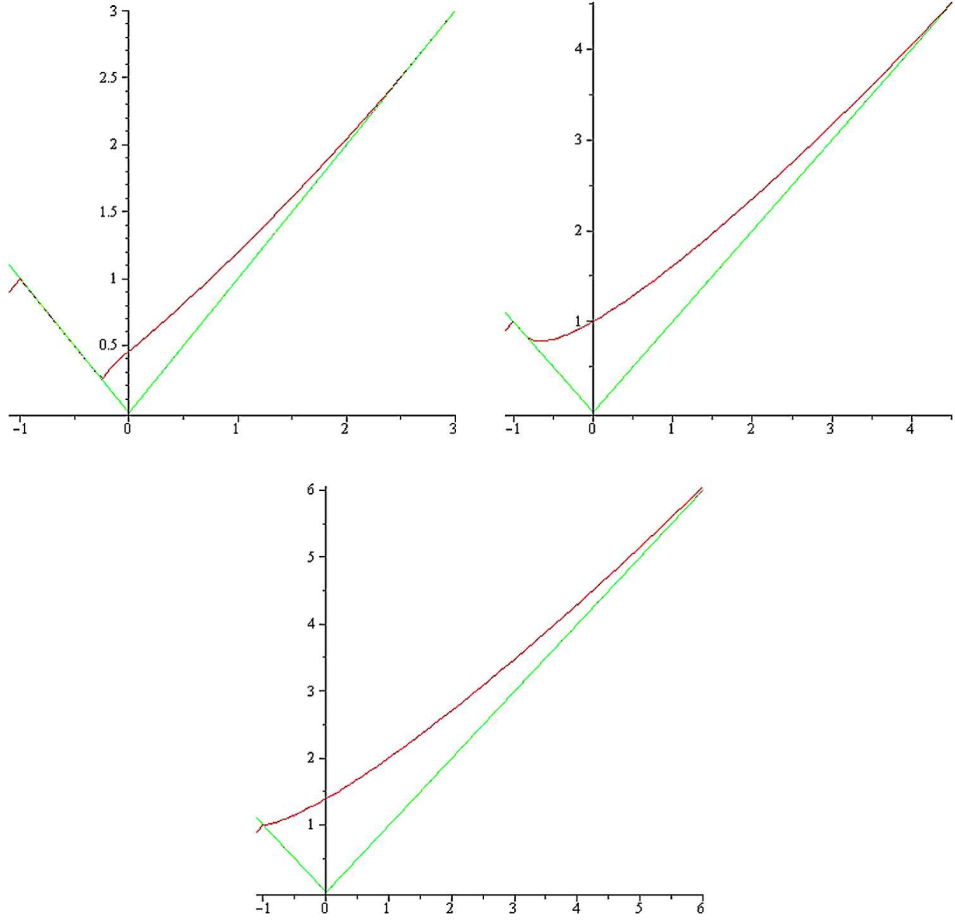}

\caption{The limit shapes for one folded $\alpha^+$-parameter with
$\alpha=1$ and $a =0.25, 1, 2$, respectively.}
\label{limitshapes1A-1}
\end{figure}

Limit shapes for $\alpha=1$ and various $a$ are shown in Figure~\ref{limitshapes1A-1}.

%
\begin{figure}[b]

\includegraphics{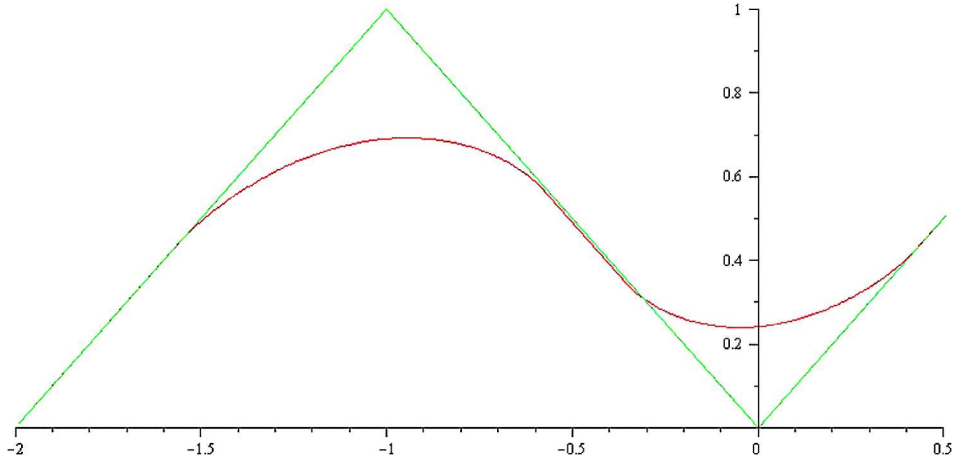}

\caption{The limit shape for two-sided Plancherel character with parameters
$\gamma_1=0.04$, $\gamma_2=0.06$.}
\label{2sPlanch}
\end{figure}

(c) \textit{One multiple $\beta^+$-parameter}.
Let us fix two positive real numbers $b$ and $\beta\le1$.
Assume that $\beta_1^+ = \beta_2^+ = \cdots= \beta_{[b N]}^+ = \beta$.

The computations in this case are equivalent to the previous one:
\begin{eqnarray*}
Q(z)&=& 1+ z(1+z) \frac{b \beta}{1+\beta z},
\\
\mathrm{Stil}^{\mathrm{multi}\mbox{-}\beta} (z) &=& \log\frac{z (1 - 2\beta) + \beta(b-1) - \sqrt
{(z-\beta(b-1))^2 + 4 b \beta^2 - 4 b \beta}}{2 \beta b - 2 \beta z}.
\end{eqnarray*}

Inside the interval $[\beta(b-1) - 2 \sqrt{b \beta(1- \beta)}; \beta
(b-1) + 2 \sqrt{b \beta(1- \beta)} ]$, the density has the following form:
\[
d^{\mathrm{multi}\mbox{-}\beta} (x) = \frac{1}{\pi} \arccos\frac{(1-2 \beta)x + \beta
(b-1)}{2 \sqrt{\beta(1-\beta) (1+z)(b-z)}}.
\]

Furthermore, as in the previous case, for some parameters $\beta$ and
$b$ there exist intervals with constant density which is equal to 1.

(d) \textit{Two-sided Plancherel character}.
Assume that $\gamma^+ = \gamma_1 N$, $\gamma^- = \gamma_2 N$ for fixed
$\gamma_1$ and $\gamma_2$, and all other parameters are equal to 0. Then
\[
Q(z) = 1 + z(1+z) \biggl( \gamma_1 - \frac{\gamma_2}{(z+1)^2} \biggr).
\]
Hence, to obtain an explicit formula for the answer one need to solve
the cubic equation
\[
\frac{z(z+1)}{(z+1) + \gamma_1 z (z+1)^2 - \gamma_2 z} = y.
\]
If $z=v_0(y)$ is the formal power series that satisfies this equation,
then the Stieltjes transform is equal to
\[
\mathrm{Stil}^{2sP} (z) = \log \biggl(1 + v_0 \biggl(
\frac{1}{z} \biggr) \biggr).
\]

An example of such a shape is shown in Figure~\ref{2sPlanch}.

In greater detail, this case was studied in \cite{BorKuan}.

(e) \textit{The case of one multiple $\alpha^+$-parameter and one multiple
$\alpha^-$-parameter}.
Assume that $\alpha_1^+ = \alpha_2^+ = \cdots= \alpha_{[a N]}^+ = \alpha
^+$ and $\alpha_1^- = \alpha_2^- = \cdots= \alpha_{[a N]}^- = \alpha
^-$. Then
\[
P'(z) = \frac{a \alpha}{1- \alpha_1 z} - \frac{ \tilde a \tilde\alpha
}{(1+z) (1+ (\tilde\alpha+1) z)}.
\]

In this case, the generating function of moments is determined by the
solution of the cubic equation
\[
\frac{z}{1 + z (z+1)  ( {a \alpha}/{(1- \alpha_1 z)} - ({
\tilde a \tilde\alpha}/{((1+z) (1+ (\tilde\alpha+1) z))} ) ) }= y.
\]

(g) \textit{The case of the continuous limit measure}.
Assume that $\alpha_i^+ = i/N$ for $i=1, \ldots, N$. It is easy to see that
\[
P'(z) = \frac{-z - \log(1-z)}{z^2}.
\]

Then the generating function of moments is determined by the solution
of the equation
\[
\frac{z^2}{-z^2 - (1+z) \log(1-z)} = y.
\]
\end{appendix}

\section*{Acknowledgements}
The authors would like to thank Richard P. Stanley for a useful remark.
The authors would like to thank the referee for valuable comments.






\printaddresses
\end{document}